\documentclass[11pt]{article}
\usepackage{latexsym,amsmath,amssymb,amsfonts}
\usepackage{diagrams,epsfig,graphics,color}
\usepackage{hyperref}
\usepackage{rotating}

\RequirePackage[mathcal,mathbf]{euler}

\usepackage{amsthm}
\newtheorem{Proposition}{Proposition}
\newtheorem{Lemma}[Proposition]{Lemma}
\newtheorem{Theorem}[Proposition]{Theorem}
\newtheorem{Corollary}[Proposition]{Corollary}

\theoremstyle{definition}
\newtheorem*{Definition}{Definition}
\newtheorem*{Notation}{Notation}

\newenvironment{Proof}{\begin{proof}}{\end{proof}}

\DeclareFontFamily{U}{rsf}{}
\DeclareFontShape{U}{rsf}{m}{n}{
  <5> <6> rsfs5 <7> <8> <9> rsfs7 <10-> rsfs10}{}
\DeclareMathAlphabet{\mathscr}{U}{rsf}{m}{n}

\renewcommand{\mathcal}{\mathscr}

\makeatletter
\def\operator@font{\sf}
\makeatother

\newcommand{\gVar}{{\mathcal{V}\textsf{ar}}}
\newcommand{\gAlg}{{\mathcal{A}\textsf{lg}}}
\newcommand{\gCat}{{\mathcal{C}\textsf{at}}}
\newcommand{\gC}{\mathcal{C}}
\newcommand{\C}{\mathbf{C}}
\newcommand{\D}{\mathbf{D}}
\newcommand{\cE}{{\mathcal E}}
\newcommand{\cF}{{\mathcal F}}
\newcommand{\cG}{{\mathcal G}}
\newcommand{\cH}{{\mathcal H}}
\newcommand{\cO}{{\mathcal O}}

\renewcommand{\k}{\textsf{k}}
\DeclareMathOperator{\Tr}{Tr}

\DeclareMathOperator{\Hom}{Hom}

\DeclareMathOperator{\Td}{td}
\DeclareMathOperator{\ch}{ch}

\DeclareMathOperator{\Ext}{Ext}
\DeclareMathOperator{\sHom}{\underline{Hom}}
\DeclareMathOperator{\bHom}{\textbf{\textsf{Hom}}}
\DeclareMathOperator{\THom}{{2\mbox{-}Hom}}
\DeclareMathOperator{\OHom}{1\mbox{-}Hom}
\DeclareMathOperator{\HH}{HH}
\DeclareMathOperator{\K}{K}

\renewcommand{\H}{\mathsf H}
\newcommand{\id}{\mathsf {Id}}
\newcommand{\Id}{\mathsf {Id}}

\newcommand{\iso}{\cong}
\newcommand{\dual}{\,{\stackrel{\chk}{\longleftrightarrow}}\,}

\newcommand{\chk}{{\scriptscriptstyle\vee}}

\newcommand{\pt}{{\mathsf{pt}}}
\newcommand{\bone}{{\mathbf 1}}
\newcommand{\adjoint}{\dashv}
\newcommand{\ra}{\rightarrow}

\newcommand{\lra}{\longrightarrow}

\newcommand{\scdot}{{\,\cdot\,}}
\newcommand{\MP}[2]{{\ensuremath{\left\langle #1,\, #2\right\rangle_{\rm M}}}}
\newcommand{\MPbig}[2]{{\ensuremath{\Biggl\langle #1,\, #2\Biggr\rangle_{\rm M}}}}
\newcommand{\SP}[2]{{\ensuremath{\left\langle #1,\, #2\right\rangle_{\rm S}}}}
\newcommand{\SPbig}[2]{{\ensuremath{\Biggl\langle #1,\, #2\Biggr\rangle_{\rm S}}}}
\renewcommand{\phi}{\varphi}
\newcommand{\SerF}{\mathsf{S}}
\newcommand{\SerK}{\Sigma}
\newcommand{\kerto}{\to}
\newcommand{\open}{\textsf{O}}
\newcommand{\closed}{\textsf{C}}
\newcommand{\TwoCob}{2\textsf{Cob}_{\text{oc}}}
\newcommand{\TwoCobL}{\TwoCob^\Lambda}

\newcommand{\circv}{\circ_{\text{v}}}
\newcommand{\circh}{\circ_{\text{h}}}
\newcommand{\one}{\mathbf{1}}
\newcommand{\OO}{\cO}
\newcommand{\blob}{\bullet}
\newcommand{\figdir}{pstex}

\newcommand{\pic}{\raisebox{-.45\height}{\input{\figdir/.pstex_t}}}

\newarrow{TeXto}----{->}
\allowdisplaybreaks

\newenvironment{Description}{\begin{list}{}{
               \setlength{\labelwidth}{2.5em}
               \setlength{\labelsep}{0.5em}
               \setlength{\itemindent}{2em}
               \setlength{\itemsep}{0.5em}
               \setlength{\rightmargin}{0em}
               \setlength{\leftmargin}{1em}}}
               {\end{list}}

\newcommand{\pd}[1]{#1^{\dagger}}
\newcommand{\rote}{{\begin{turn}{270}$\!\!\!\epsilon\ $\end{turn}}}

\author{%
Andrei C\u ald\u araru\thanks{Mathematics Department,
University of Wisconsin--Madison, 480 Lincoln Drive, Madison, WI
53706--1388, USA, {\em e-mail: }{\tt andreic@math.wisc.edu}}\ \  and Simon
Willerton\thanks{Department of Pure Mathematics, University of Sheffield,
Hicks Building, Hounsfield Road, Sheffield, S3 7RH, UK, {\em e-mail:
}{\tt S.Willerton@sheffield.ac.uk}}}

\title{The Mukai pairing, I: a categorical approach}

\date{July 12, 2007}

\begin{document}

\maketitle

\begin{abstract}
  We study the Hochschild homology of smooth spaces, emphasizing the
  importance of a pairing which generalizes Mukai's pairing on the
  cohomology of K3 surfaces.  We show that integral transforms between
  derived categories of spaces induce, functorially, linear maps on
  homology.  Adjoint functors induce adjoint linear maps with respect
  to the Mukai pairing.  We define a Chern character with values in
  Hochschild homology, and we discuss analogues of the
  Hirzebruch-Riemann-Roch theorem and the Cardy Condition from
  physics.  This is done in the context of a 2-category which has
  spaces as its objects and integral kernels as its 1-morphisms.
\end{abstract}

\section*{Introduction}
\label{sec:intro}

\addcontentsline{toc}{section}{Introduction}

The purpose of the present paper is to introduce the Mukai pairing on
the Hochschild homology of smooth, proper spaces.  This pairing is the
natural analogue, in the context of Hochschild theory, of the
Poincar\'e pairing on the singular cohomology of smooth manifolds.

Our approach is categorical.  We start with a geometric category,
whose objects will be called {\em spaces}.  For a space $X$ we define
its Hochschild homology which is a graded vector space $\HH_\blob(X)$
equipped with the non-degenerate \emph{Mukai pairing}.  We show that
this structure satisfies a number of properties, the most important of
which are {\em functoriality} and {\em adjointness}.

The advantage of the categorical approach is that the techniques we
develop apply in a wide variety of geometric situations, as long as an
analogue of Serre duality is satisfied.  Examples of categories for
which our results apply include compact complex manifolds, proper
smooth algebraic varieties, proper Deligne-Mumford stacks for which
Serre duality holds, representations of a fixed finite group, and
compact ``twisted spaces'' in the sense of~\cite{Cal}.  We expect the
same construction to work for categories of Landau-Ginzburg
models~\cite{Orl}, but at the moment it is not known whether Serre
duality holds for these.

\subsection*{The Hochschild structure}
\label{subsec:hochstr}
In order to define the Hochschild structure of a space we need
notation for certain special kernels which will play a fundamental
role in what follows.  For a space $X$, denote by $\Id_X$ and
$\SerK_X^{-1}$ the objects of $\D(X\times X)$ given by
\[ \Id_X := \Delta_* \cO_X \quad\mbox{and}\quad \SerK_X^{-1} :=
\Delta_* \omega_X^{-1}[-\dim X], \] 
where $\Delta:X \ra X\times X$ is the diagonal map, and
$\omega_X^{-1}$ is the anti-canonical line bundle of $X$.  When
regarded as kernels, these objects induce the identity functor and the
inverse of the Serre functor on $\D(X)$, respectively.  We shall see
in the sequel that $\Id_X$ can be regarded as the identity 1-morphism
of $X$ in a certain 2-category $\gVar$.

The {\em Hochschild structure} of the space $X$ then consists of the
following data:
\begin{itemize}
\item the graded ring $\HH^\blob(X)$, the Hochschild cohomology ring of
  $X$, whose $i$-th graded piece is defined as
  \[ \HH^i(X) := \Hom_{\D(X\times X)}^i(\Id_X, \Id_X); \]
\item the graded left $\HH^\blob(X)$-module $\HH_\blob(X)$, the
  Hochschild homology module of $X$, defined as
  \[ \HH_i(X) := \Hom_{\D(X\times X)}^{-i}(\SerK_X^{-1}, \Id_X); \]
\item a non-degenerate graded pairing $\MP{{-}}{{-}}$ on
  $\HH_\blob(X)$, the {\em generalized Mukai pairing}.
\end{itemize}

The above definitions of Hochschild homology and cohomology agree with
the usual ones for quasi-projective schemes (see~\cite{CalHH2}).  The
pairing is named after Mukai, who was the first to introduce a pairing
satisfying the main properties below, on the total cohomology of
complex K3 surfaces~\cite{MukK3}.

\subsection*{Properties of the Mukai pairing}
The actual definition of the Mukai pairing is quite complicated and is
given in Section~\ref{sec:adjoint}.  We can, however, extricate the
fundamental properties of Hochschild homology and of the Mukai
pairing.
\begin{Description}
\item[Functoriality] Integral kernels induce, in a functorial way,
  linear maps on Hochschild homology. Explicitly, to any integral
  kernel \( \Phi\in \D(X\times Y) \) we associate, in
  Section~\ref{def:defmaps}, a linear map of graded vector spaces
  \[ \Phi_* \colon \HH_\blob(X) \ra \HH_\blob(Y), \]
  and this association is functorial with respect to composition of
  integral kernels (Theorem~\ref{thm:functorial}).
\item[{Adjointness}] For any adjoint pair of integral kernels
  $\Psi\adjoint \Phi$, the induced maps on homology are themselves
  adjoint with respect to the Mukai pairing:
  \[ \MP{\Psi_* v}{w} = \MP{v}{\Phi_* w}\]
  for $v\in \HH_\blob(Y)$, $w \in \HH_\blob(X)$
  (Theorem~\ref{thm:adjointness}).
\end{Description}
The following are then consequences of the above basic properties:
\begin{Description}
\item[{Chern character}] In all geometric situations there is a
  naturally defined object $\bone \in \HH_0(\pt)$.  An element $\cE$
  in $\D(X)$ can be thought of as the kernel of an integral transform
  $\pt \ra X$, and using functoriality of homology we define a Chern
  character map
  \[ \ch\colon\K_0(X) \ra \HH_0(X), \quad \ch(\cE) = \cE_*(\bone).  \]
  For a smooth proper variety the Hochschild-Kostant-Rosenberg
  isomorphism identifies $\HH_0(X)$ and $\bigoplus_p \H^{p,p}(X)$; our
  definition of the Chern character matches the usual one under this
  identification~\cite{CalHH2}.
\item[Semi-Hirzebruch-Riemann-Roch Theorem] For $\cE, \cF \in \D(X)$
  we have
  \[ \MP{\ch(\cE)}{\ch(\cF)} = \chi(\cE, \cF) = \sum_i (-1)^i \dim
  \Ext_X^i(\cE, \cF). \]
\item[Cardy Condition] The Hochschild structure appears naturally in
  the context of open-closed topological quantum field theories
  (TQFTs).  The Riemann-Roch theorem above is a particular case of a
  standard constraint in these theories, the Cardy Condition.  We
  briefly discuss open-closed TQFTs, and we argue that the natural
  statement of the Cardy Condition in the B-model open-closed TQFT is
  always satisfied, even for spaces which are not Calabi-Yau
  (Theorem~\ref{thm:cardy}).
\end{Description}

\subsection*{The 2-categorical perspective}
In order to describe the functoriality of Hochschild homology it is
useful to take a macroscopic point of view using a 2-category called
$\gVar$.  One way to think of this 2-category is as something half-way
between the usual category consisting of spaces and maps, and $\gCat$,
the 2-category of (derived) categories, functors and natural
transformations.  The 2-category $\gVar$ has spaces as its objects,
has objects of the derived category $\D(X\times Y)$ --- considered as
integral kernels --- as its 1-morphisms from $X$ to $Y$, and has
morphisms in the derived category as its
2-morphisms.

One consequence of thinking of spaces in this 2-category is that
whereas in the usual category of spaces and maps two spaces are
equivalent if they are isomorphic, in $\gVar$ two spaces are
equivalent precisely when they are Fourier-Mukai partners.  This is
the correct notion of equivalence in many circumstances, thus making
$\gVar$ an appropriate context in which to work.

This point of view is analogous to the situation in Morita theory in
which the appropriate place to work is not the category of algebras
and algebra morphisms, but rather the 2-category of algebras,
bimodules and bimodule morphisms.  In this 2-category two algebras are
equivalent precisely when they are Morita equivalent, which again is
the pertinent notion of equivalence in many situations.

Many facts about integral transforms can be stated very elegantly as
facts about the 2-category $\gVar$.  For example, the fact that every
integral transform between derived categories has both a left and
right adjoint is an immediate consequence of the more precise fact ---
proved exactly the same way --- that every integral kernel has both a
left and right adjoint in $\gVar$.  Here the definition of an adjoint
pair of 1-morphisms in a 2-category is obtained from one of the
standard definitions of an adjoint pair of functors by everywhere replacing the
word `functor' by the word `1-morphism' and the words `natural
transformation' by the word `2-morphism'.

The Hochschild cohomology of a space $X$ has a very natural
description in terms of the 2-category $\gVar$: it is the ``second
homotopy group of $\gVar$ based at $X$'', which means that it is
$\THom_\gVar(\Id_X,\Id_X)$, the set of 2-morphisms from the identity
1-morphism at $X$ to itself.  Unpacking this definition for $\gVar$
one obtains precisely $\Ext^\blob_{X\times X}(\cO_\Delta,\cO_\Delta)$,
one of the standard definitions of Hochschild cohomology.  By analogy
with homotopy groups, given a kernel $\Phi\colon X\to Y$, i.e., a
``path'' in $\gVar$, one might expect an induced map $\HH^\blob(X)\to
\HH^\blob(Y)$ obtained by ``conjugating with $\Phi$''.  However, this
does not work, as the analogue of the ``inverse path to $\Phi$''
needed is a simultaneous left and right adjoint of $\Phi$, and such a
thing does not exist in general as the left and right adjoints of
$\Phi$ differ by the \emph{Serre kernels} of $X$ and $Y$.

The Hochschild homology $\HH_\blob(X)$ of a space $X$ can be given a
similar natural definition in terms of $\gVar$ --- it is
$\THom_\gVar(\SerK^{-1}_X,\Id_X)$ the set of 2-morphisms from the
inverse Serre kernel of $X$ to the identity 1-morphism at $X$.
In this case, the idea of ``conjugating by a kernel $\Phi\colon X \to
Y$'' does work as the Serre kernel in the definition exactly
compensates the discrepancy between the left and right adjoints of
$\Phi$.

The functoriality of Hochschild homology can be expressed by saying
that $\HH_\blob$ is a functor into the category of vector spaces from
the Grothen\-dieck category of the 2-category $\gVar$ (i.e., the
analogue of the Grothendieck group of a 1-category) whose objects are
spaces and whose morphisms are \emph{isomorphism classes} of kernels.
One aspect of this which we do not examine here is related to the fact
that this Grothendieck category is actually a monoidal category with
certain kinds of duals for objects and morphisms, and that Hochschild
homology is a monoidal functor.  The Mukai pairing is then a
manifestation of the fact that spaces are self-dual in this
Grothendieck category.  Details will have to appear elsewhere.

There is an alternative categorical approach to defining Hoch\-schild
homology and cohomology.  This approach uses the notion of enhanced
triangulated categories of Bondal and Kapranov~\cite{BonKap}, which
are triangulated categories arising as homotopy categories of
differential-graded (dg) categories.  In~\cite{Toe}, To\"en argued that
the Hochschild cohomology $\HH^\blob(X)$ of a space $X$ can be
regarded as the cohomology of the dg-algebra of dg-natural
transformations of the identity functor on the dg-enhancement of
$\D(X)$.  It seems reasonable to expect that a similar construction
can be used to define the Hochschild homology $\HH_\blob(X)$ as
dg-natural transformations from the inverse of the Serre functor to
the identity.  However, since the theory of Serre functors for
dg-categories is not yet fully developed, we chose to use the language
of the 2-category $\gVar$, where all our results can be made precise.

\subsection*{String diagram notation}
As 2-categories are fundamental to the functoriality, and they are
fundamentally 2-dimensional creatures, we adopt a 2-dimensional
notation.  The most apt notation in this situation appears to be that
of \emph{string diagrams}, which generalizes the standard notation
used for monoidal categories in quantum topology.  String diagrams are
Poincar\'e dual to the usual arrow diagrams for 2-categories.  The
reader unfamiliar with these ideas should be aware that the pictures
scattered through this paper form rigorous notation and are not just
mnemonics.

\subsection*{Note} This paper supersedes the unpublished paper
\cite{CalMPI}, in which it was stated that hopefully the correct
categorical context could be found for the results therein.  This
paper is supposed to provide the appropriate context.

\subsection*{Synopsis} The paper is structured as follows.  The first
section is devoted to the study of integral transforms and of the
2-category $\gVar$.  In the next section we review the Serre functor and
Serre trace on the derived category, and we use these in
Section~\ref{sec:adjoints} to study adjoint kernels in $\gVar$.  In
Section~\ref{sec:defmaps} we introduce the maps between Hochschild
homology groups associated to a kernel, and we examine their
functorial properties.  The Mukai pairing is defined in
Section~\ref{sec:adjoint}, where we also prove its compatibility with
adjoint functors.  In Section~\ref{sec:chernchar} we define the Chern
character and we prove the Semi-Hirzebruch-Riemann-Roch theorem.  We conclude with
Section~\ref{sec:cardy} where we review open-closed TQFTs, and we
discuss the Cardy Condition.  An appendix contains some of the more
technical proofs.

\subsection*{Notation} Throughout this paper $\k$ will denote an
algebraically closed field of characteristic zero and $\D(X)$ will
denote the bounded derived category of coherent sheaves on $X$.
Categories will be denoted by bold letters, such as $\C$, and the
names of 2-categories will have a script initial letter, such as
$\gVar$.

 \subsection*{The base category of spaces}
We fix for the remainder of the paper a geometric category, whose
objects we shall call \emph{spaces}.  It is beyond the purpose of this
paper to list the axioms that this category needs to satisfy, but the
following categories can be used:
\begin{itemize}
\item smooth projective schemes over $\k$;
\item smooth proper Deligne-Mumford stacks over $\k$;
\item smooth projective schemes over $\k$, with an action of a fixed
  finite group $G$, along with $G$-equivariant morphisms;
\item twisted spaces in the sense of~\cite{Cal}, i.e., smooth
  projective schemes over $\k$, enriched with a sheaf of Azumaya
  algebras.
\end{itemize}
For any space $X$ as above, the category of coherent sheaves on $X$
makes sense, and the standard functors (push-forward, pull-back,
sheaf-hom, {etc}.)  are defined and satisfy the usual compatibility
relations.

\subsection*{Acknowledgments} We have greatly benefited from
conversations with Jonathan Block, Tom Bridgeland, Andrew Kresch,
Aaron Lauda, Eyal Markman, Mircea Musta\c t\u a, Tony Pantev, Justin
Roberts, Justin Sawon, and Sarah Witherspoon. Many of the ideas in
this work were inspired by an effort to decipher the little known but
excellent work~\cite{Mar} of Nikita Markarian. Greg Moore and Tom
Bridgeland suggested the connection between the Semi-Hirzebruch-Riemann-Roch
theorem and the Cardy Condition in physics.

AC's initial work on this project was supported by an NSF
post-doctoral fellowship, and by travel grants and hospitality from
the University of Pennsylvania, the University of Salamanca, Spain,
and the Newton Institute in Cambridge, England.  AC's current work is
supported by the National Science Foundation under Grant No.\
DMS-0556042.  SW has been supported by a WUN travel bursary and a
Royal Society Conference grant.

\section{The 2-category of kernels}

In this section we introduce the 2-category $\gVar$, which provides
the natural context for the study of the structure of integral
transforms between derived categories of spaces.  The objects of
$\gVar$ are spaces, 1-morphisms are kernels of integral transforms,
and 2-morphisms are maps between these kernels.  Before introducing
$\gVar$ we remind the reader of the notion of a 2-category and we
explain the string diagram notation of which we will have much use.

 \subsection{A reminder on 2-categories.}
We will review the notion of a 2-category at the same time as
introducing the notation we will be using. Recall that a 2-category
consists of three levels of structure: objects; 1-morphisms between
objects; and 2-morphisms between 1-morphisms.  It is worth mentioning
a few examples to bear in mind during the following exposition.

\begin{itemize}
\item[1.] The first example is the 2-category $\gCat$ of categories,
  functors and natural transformations.
\item[2.] The second example is rather a family of examples.  There is
  a correspondence between 2-categories with one object $\star$ and
  monoidal categories.  For any monoidal category the objects and
  morphisms give respectively the 1-morphisms and 2-morphisms of the
  corresponding 2-category.
\item[3.] The third example is the 2-category $\gAlg$ with algebras
  over some fixed commutative ring as its objects, with the set of
  $A$-$B$-bimodules as its 1-morphisms from $A$ to $B$, where
  composition is given by tensoring over the intermediate algebra, and
  with bimodule morphisms as its 2-morphisms.
\end{itemize}

There are various ways of notating 2-categories: the most common way
is to use \emph{arrow diagrams}, however the most convenient way for
the ideas in this paper is via \emph{string diagrams} which are
Poincar\'e dual to the arrow diagrams.  In this subsection we will
draw arrow diagrams on the left and string diagrams on the right to
aid the reader in the use of string diagrams.

Recall the idea of a 2-category.  For any pair of objects $X$ and $Y$
there is a collection of morphisms $\OHom(X,Y)$; if
$\Phi\in\OHom(X,Y)$ is a 1-morphism then it is drawn as below.
\[\pic{AD1}\qquad\pic{SD1}\]
These 1-dimensional pictures will only appear as the source and target
of 2-morphism, i.e., the top and bottom of the 2-dimensional pictures
we will be using.  In general 1-morphisms will be denoted by their
identity 2-morphisms, see below.

If $\Phi,\Phi'\in\OHom(X,Y)$ are parallel 1-morphisms --- meaning
simply that they have the same source and target --- then there is a
set of 2-morphisms $\THom(\Phi,\Phi')$ from $\Phi$ to $\Phi'$.  If
$\alpha\in\THom(\Phi,\Phi')$ is a 2-morphism then it is drawn as
below.
\[\pic{AD2}\qquad\pic{SD2}\]
At this point make the very important observation that diagrams are read from right to
left and from bottom to top.

There is a \emph{vertical composition} of 2-morphisms so that if
$\alpha\colon\Phi\Rightarrow\Phi'$ and
$\alpha'\colon\Phi'\Rightarrow\Phi''$ are 2-morphisms then the
vertical composite $\alpha'\circv\alpha\colon\Phi\Rightarrow\Phi''$ is
defined and is denoted as below.
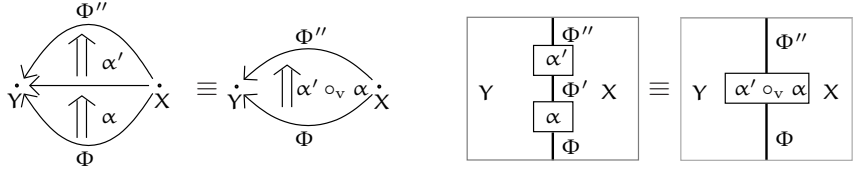
\[\pic{AD3a}\equiv\pic{AD3b}\qquad\pic{SD3}\equiv\pic{SD3b}\]
This vertical composition is strictly associative so that
$(\alpha''\circv\alpha')\circv\alpha
=\alpha''\circv(\alpha'\circv\alpha)$ whenever the three 2-morphisms
are composable.  Moreover, there is an identity 2-morphism
$\id_\Phi\colon \Phi\Rightarrow \Phi$ for every 1-morphism $\Phi$ so
that $\alpha\circv\id_\Phi=\alpha=\id_{\Phi'}\circv\alpha$ for every
2-morphism $\alpha\colon\Phi\Rightarrow\Phi'$.  This means that for
every pair of objects $X$ and $Y$, the 1-morphisms between them are the
objects of a category $\bHom(X,Y)$, with the 2-morphisms forming the
morphisms.  In the string notation the identity 2-morphisms are usually just
drawn as straight lines.
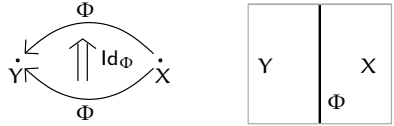
\[\pic{AD4}\qquad\pic{SD4}\]

There is also a composition for 1-morphisms, so if $\Phi\colon X\to Y$
and $\Psi \colon Y\to Z$ are 1-morphisms then the composite $\Psi\circ
\Phi\colon X\to Z$ is defined and is denoted as below.
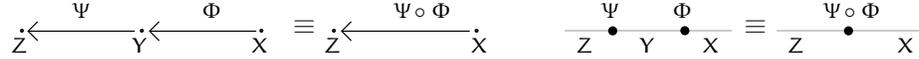
\[\pic{AD5a}\equiv\pic{AD5b}\qquad\pic{SD5}\equiv\pic{SD5b}\]
Again, these pictures will only appear at the top and bottom of
2-morphisms.  This composition of 1-morphisms is not required to be
strictly associative, but it is required to be associative up to a
\emph{coherent} 2-isomorphism.  This means that for every composable
triple $\Theta$, $\Psi$ and $\Phi$ of 1-morphisms there is a specified
2-isomorphism $(\Theta\circ\Psi)\circ\Phi
\stackrel\sim\Longrightarrow\Theta\circ(\Psi\circ\Phi)$ and these
2-isomorphisms have to satisfy the so-called pentagon coherency
condition which ensures that although $\Theta\circ\Psi\circ\Phi$ is
ambiguous, it can be taken to mean either $(\Theta\circ\Psi)\circ\Phi$
or $\Theta\circ(\Psi\circ\Phi)$ without confusion.  The up-shot of
this is that parentheses are unnecessary in the notation.

Each object $X$ also comes with an identity 1-morphism $\id_X$, but
again, in general, one does not have equality of ${\id_Y}\circ \Phi$,
$\Phi$ and $\Phi\circ\id_X$, but rather the identity 1-morphisms come
with coherent 2-isomorphisms ${\id_Y}\circ\Phi\stackrel\sim
\Rightarrow\Phi$, and $\Phi\circ\id_X\stackrel\sim\Rightarrow\Phi$.
Again this means that in practice the identities can be neglected in
the notation: so although we could denote the identity 1-morphism
with, say, a dotted line, we choose not to. This is illustrated below.
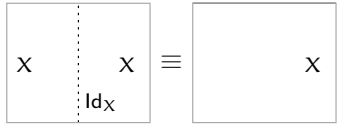
\[\pic{SD6a}\equiv\pic{SD6b}\]

A \emph{strict} 2-category is one in which the coherency
2-isomorphisms for associativity and identities are themselves all
identities.  So the 2-category $\gCat$ of categories, functors and
natural transformations is a strict 2-category.

The last piece of structure that a 2-category has is the
\emph{horizontal composition} of 2-morphisms.  If $\Phi,\Phi'\colon
X\to Y$ and $\Psi,\Psi'\colon Y\to Z$ are 1-morphisms, and
$\alpha\colon \Phi\Rightarrow\Phi'$ and $\beta\colon
\Psi\Rightarrow\Psi'$ are 2-morphisms, then $\beta\circh \alpha\colon
\Psi\circ\Phi\Rightarrow\Psi'\circ\Phi'$ is defined, and is notated as
below.
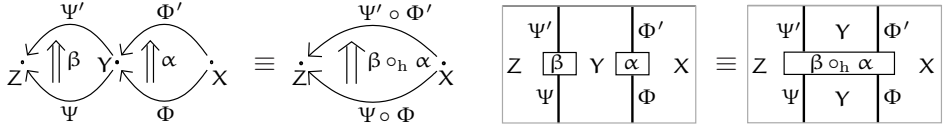
\[\pic{AD7a}\equiv\pic{AD7b}\quad\pic{SD7a}\equiv\pic{SD7b}\]
The horizontal and vertical composition are required to obey the
\emph{interchange law}.
\[(\beta'\circv\beta)\circh(\alpha'\circv\alpha)=
(\beta'\circh\alpha')\circv(\beta\circh\alpha).\]
This means that the following diagrams are unambiguous.
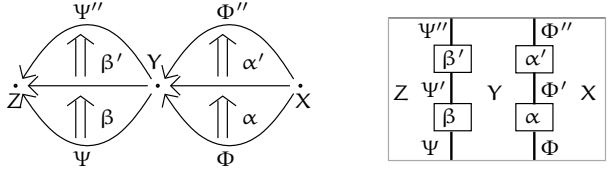
\[\pic{AD8}\qquad\pic{SD8}\]
It also means that 2-morphisms can be `slid past' each other in the
following sense.
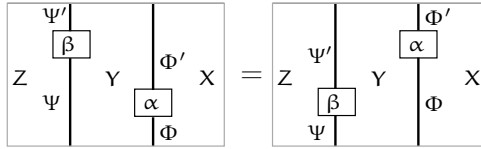
\[\pic{SD9a}=\pic{SD9b}\]

From now on, string diagrams will be drawn without the grey borders, and labels will be omitted if they are clear from the context.

 \subsection{The 2-category $\gVar$.}
The 2-category $\gVar$, of spaces and integral kernels, is defined as
follows.  The objects are spaces, as defined in the introduction, and the hom-category $\bHom_\gVar(X,Y)$ from a space $X$
to a space $Y$ is the derived category $\D(X\times Y)$, which is to be
thought of as the category of integral kernels from $X$ to $Y$.
Explicitly, this means that the 1-morphisms in $\gVar$ from $X$ to $Y$
are objects of $\D(X\times Y)$ and the 2-morphisms from $\Phi$ to
$\Phi'$ are morphisms in $\Hom_{\D(X\times Y)}(\Phi,\Phi')$, with
vertical composition of 2-morphisms just being usual composition in
the derived category.

Composition of 1-morphisms in $\gVar$ is defined using the convolution
of integral kernels: if $\Phi\in \D(X\times Y)$ and $\Psi\in
\D(Y\times Z)$ are 1-morphisms then define the convolution
$\Psi\circ\Phi\in\D(X\times Z)$ by
\[ \Phi\circ \Psi := \pi^{}_{XZ, *} ( \pi_{YZ}^* \Psi \otimes
\pi_{XY}^* \Phi), \] where $\pi_{XZ}$, $\pi_{XY}$ and $\pi_{YZ}$ are the
projections from $X\times Y\times Z$ to the appropriate factors.  The
horizontal composition of 2-morphisms is similarly defined.  Finally,
the identity 1-morphism $\id_X\colon X\ra X$ is given by
$\cO_{\Delta}\in \D(X\times X)$, the structure sheaf of the diagonal
in $X\times X$.

The above 2-category is really what is at the heart of the study of
integral transforms, and it is entirely analogous to $\gAlg$, the
2-category of algebras described above.  For example, the Hochschild
cohomology groups of a space $X$ arise as the second homotopy groups
of the 2-category $\gVar$, at $X$:
\begin{align*}
  \HH^\blob(X) &:=\Ext_{X\times X}^\blob(\cO_\Delta,\cO_\Delta)
  \cong \Hom^\blob_{D(X\times X)}(\cO_\Delta,\cO_\Delta)\\
  &=:\THom_{\gVar}(\Id_X,\Id_X).
\end{align*}

There is a 2-functor from $\gVar$ to $\gCat$ which encodes integral
transforms: this 2-functor sends each space $X$ to its derived
category $\D(X)$, sends each kernel $\Phi\colon X\to Y$ to the
corresponding integral transform $\overline \Phi\colon D(X)\to D(Y)$,
and sends each map of kernels to the appropriate natural
transformation.  Many of the statements about integral transforms have
better formulations in the language of the 2-category $\gVar$.

\section{Serre functors}
\label{sec:serre}

In this section we review the notion of the Serre functor on $\D(X)$
and then show how to realise the Serre functor on the derived category
$\D(X\times Y)$ using 2-categorical language.

 \subsection{The Serre functor on $\D(X)$.}
\label{subsec:Serretrace}
If $X$ is a space then we consider the functor
\[\SerF\colon\D(X)\to \D(X);\qquad \cE\mapsto \omega_X[\dim
X]\otimes \cE,\] 
where $\omega_X$ is the canonical line bundle of $X$.  Serre duality
then gives natural, bifunctorial isomorphisms
\[ \eta_{\cE,\cF}\colon\Hom_{\D(X)}(\cE, \cF) \stackrel{\sim}{\lra}
\Hom_{\D(X)}(\cF, \SerF \cE)^\chk \]
for any objects $\cE,\cF\in \D(X)$, where ${-}^\chk$ denotes the dual
vector space.

A functor such as $\SerF$, together with isomorphisms as above, was
called a \emph{Serre functor} by Bondal and Kapranov
\cite{BonKapSerre} (see also~\cite{ReiVan}).  From this data, for any
object $\cE\in\D(X)$, define the \emph{Serre trace} as follows:
\[ \Tr\colon\Hom(\cE, \SerF \cE) \ra \k;\qquad \Tr(\alpha) := \eta_{\cE,\cE}(\id_\cE)(\alpha). \]
Note that from this trace we can recover $\eta_{\cE,\cF}$ because
\[ \eta_{\cE,\cF}(\alpha)(\beta)=\Tr(\beta\circ\alpha). \]
We also have the commutativity identity
\[ \Tr(\beta\circ \alpha) = \Tr(\SerF \alpha \circ \beta). \]
Yet another way to encode this data is as a perfect pairing, the
\emph{Serre pairing}:
\[\SP{{-}}{{-}}\colon \Hom(\cE,\cF)\otimes \Hom(\cF,\SerF \cE)\to
\k; \qquad \SP{\alpha}{\beta}:=\Tr(\beta\circ \alpha).\]

 \subsection{Serre kernels and the Serre functor for $\D(X\times
  Y)$.}
\label{subsec:serXY}
We are interested in kernels and the 2-category $\gVar$, so are
interested in Serre functors for product spaces $X\times Y$, and
these have a lovely description in the 2-categorical language.
We can now define one of the key objects in this paper.

\begin{Definition}
For a space $X$, the \emph{Serre kernel} $\SerK_X\in \OHom_\gVar(X,X)$
is defined to be $\Delta_*\omega_X[\dim X]\in\D(X\times X)$, the
kernel inducing the Serre functor on $X$.  Similarly the
\emph{anti-Serre kernel} $\SerK^{-1}_X\in \OHom_\gVar(X,X)$ is defined
to be $\Delta_*\omega^{-1}_X[-\dim X]\in\D(X\times X)$.
\end{Definition}

\begin{Notation}
In string diagrams the Serre kernel will be denoted by a dashed-dotted
line, while the anti-Serre kernel will be denoted by a dashed-dotted
line with a horizontal bar through it (see the pictures in
Section~\ref{subsec:leftrightadjoints}).
\end{Notation}
\noindent
The Serre kernel can now be used to give a natural description, in the
2-category language, of the Serre functor on the product $X\times Y$.

\begin{Proposition}
For spaces $X$ and $Y$ the Serre functor $\SerF_{X\times Y}\colon
\D(X\times Y)\to \D(X\times Y)$
can be taken to be $\SerK_Y\circ{-}\circ\SerK_X$.
\end{Proposition}

\begin{Proof}
The Serre functor on $\D(X\times Y)$ is given by
\[\SerF_{X\times Y}(\Phi)= \Phi\otimes \pi_X^*\omega_X\otimes
\pi_Y^*\omega_Y[\dim X+\dim Y].\]
However, observe that if $\Phi\in \D(X\times Y)$ and $\cE\in \D(X)$
then
\[\Phi\circ \Delta_*\cE\cong \Phi\otimes \pi_X^* \cE,\]
where $\pi_X\colon X\times Y \to X$ is the projection.  This is just a
standard application of the base-change and projection formulas.
Similarly if $\cF\in \D(Y)$ then \( \Delta_*\cF\circ\Phi\cong \pi_Y^*
\cF\otimes\Phi.\) From this the Serre functor can be written as
\[\SerF_{X\times Y}(\Phi)= \SerK_Y\circ\Phi\circ \SerK_X.\qedhere\]
\end{Proof}

\noindent
This means that the Serre trace map on $X\times Y$ is a map
\[ \Tr\colon\THom_\gVar(\Phi, \SerK_Y \circ\Phi\circ\SerK_X) \ra \k \]
which can be pictured as
\[\Tr\left(\pic{Trace}\right) \in \k,\]
where the Serre kernel is denoted by the dashed-dotted line.

We will see below that we have `partial trace' operations which the
Serre trace factors through.

\section{Adjoint kernels}
\label{sec:adjoints}

The reader is undoubtably familiar with the notion of adjoint
functors.  It is easy and natural to generalize this from the context
of the 2-category $\gCat$ of categories, functors and natural
transformations to the context of an arbitrary 2-category.  In this
section it is shown that every kernel, considered as a 1-morphism in
the 2-category $\gVar$, has both a left and right adjoint: this is a
consequence of Serre duality, and is closely related to the familiar
fact that every integral transform functor has both a left and right
adjoint functor.

Using these notions of left and right adjoints we define \emph{left
  partial trace maps}, and similarly \emph{right partial trace
  maps}. These can be viewed as partial versions of the Serre trace
map. This construction is very much the heart of the paper.

 \subsection{Adjunctions in 2-categories.}
The notion of an adjunction in a 2-category simultaneously generalizes
the notion of an adjunction between functors and the notion of a
duality between objects of a monoidal category.  As it is the former
that arises in the context of integral transforms, we will use that as
the motivation, but will come back to the latter below.

The most familiar definition of adjoint functors is as follows.  For
categories $\C$ and $\D$, an adjunction $\Psi\dashv \Phi$ between
functors $\Psi\colon \D\to \C$ and $\Phi\colon \C\to \D$ is the
specification of a natural isomorphism
\[t_{a,b}\colon \Hom_\C(\Psi(a),b)\xrightarrow{\sim}
\Hom_\D(a,\Phi(b))\]
for any $a\in \D$ and $b\in\C$.

It is well known (see \cite[page 91]{GelfandManin}) that this definition is
equivalent to an alternative definition of adjunction which consists
of the specification of unit and counit natural morphisms, namely
\[\eta\colon \Id_\D \Rightarrow \Phi\circ \Psi
\quad\text{and}\quad \epsilon\colon \Psi\circ \Phi\Rightarrow
\Id_\C, \]
such that the composite natural transformations
\[ \Psi \stackrel{\id_\Psi\circ\eta}\Longrightarrow
\Psi\circ\Phi\circ\Psi \stackrel{\epsilon\circ\Id_\Psi}\Longrightarrow
\Psi \quad\text{and}\quad \Phi
\stackrel{\eta\circ\Id_\Phi}\Longrightarrow \Phi\circ\Psi\circ\Phi
\stackrel{\Id_\Phi\circ\epsilon}\Longrightarrow \Phi
\]
are respectively the identity natural transformation on $\Psi$ and
the identity natural transformation on $\Phi$.

It is straight forward to translate between the two different
definitions of adjunction.  Given $\eta$ and $\epsilon$ as above, define
$t_{a,b}\colon\Hom_\C(\Psi(a),b)\xrightarrow{\sim}
\Hom_\D(a,\Phi(b))$ by $t_{a,b}(f):=\Phi(f)\circ \eta_a$.  The inverse
of $t_{a,b}$ is defined similarly using the counit $\epsilon$.
Conversely, to get the unit and counit from the natural isomorphism of
hom-sets, define $\eta_a:=t_{a,\Psi(a)}(\id_{\Psi(a)})$ and define
$\epsilon$ similarly.

The definition involving the unit and counit is stated purely in terms
of functors and natural transformations --- without mentioning objects
--- thus it generalizes immediately to arbitrary 2-categories.

\begin{Definition}
If $\gC$ is a 2-category, $X$ and $Y$ are objects of $\gC$, and
\[ \Phi\colon X \ra Y\quad \text{and}\quad \Psi\colon Y \ra X \]
are 1-morphisms, then an \emph{adjunction} between $\Phi$ and $\Psi$
consists of two 2-morphisms
\[ \eta\colon \Id_Y \Rightarrow \Phi\circ\Psi \quad\text{and}\quad
\epsilon\colon \Psi\circ\Phi \Rightarrow \Id_X, \]
such that
\[(\epsilon\circh\Id_\Psi)\circv(\Id_\Psi\circh\eta)=\Id_\Psi
\quad\text{and}\quad(\Id_\Phi\circh\epsilon)\circv(\eta\circh\Id_\Phi)=\Id_\Phi.\]
Given such an adjunction we write $\Psi\dashv\Phi$.
\end{Definition}

It is worth noting that this also generalizes the notion of duality in
a monoidal category, that is to say two objects are dual in a monoidal
category if and only if the corresponding 1-morphisms are adjoint in
the corresponding 2-category-with-one-object.  Indeed, taking this
point of view, May and Sigurdsson \cite{MaySigurdsson} refer to what
is here called adjunction as duality.

It is at this point that the utility of the string diagram notation
begins to be seen.  Given an adjunction $\Psi\dashv\Phi$ the counit
$\epsilon \colon \Psi\circ \Phi \Rightarrow\Id$ and the unit
$\eta\colon \Id\Rightarrow\Phi\circ \Psi$ can be denoted as follows:
\[\pic{NaiveCounit2Adjunction};\qquad \text{and}\qquad
\pic{NaiveUnit2Adjunction}.\]
However, adopting the convention of denoting the identity one-morphism
by omission, it is useful just to draw the unit and counit as a cup
and a cap respectively:
\[\pic{Counit2Adjunction}:=\pic{NaiveCounit2Adjunction};\qquad
\text{and}\qquad \pic{Unit2Adjunction}:= \pic{NaiveUnit2Adjunction}.\]
The relations become the satisfying
\[\pic{SnakePsi}=\pic{IdPsi}\quad\text{and}\quad\pic{SnakePhi}=\pic{IdPhi}.\]

Adjunctions in 2-categories, as defined above, do correspond to
isomorphisms of certain hom-sets but in a different way to the
classical notion of adjunction.  Namely, if $\Theta\colon Z\to Y$ and
$\Xi\colon Z\to X$ are two other 1-morphisms, then an adjunction
$\Psi\dashv \Phi$ as above gives an isomorphism
\begin{align*}
  \THom(\Psi\circ\Theta,\Xi)&\xrightarrow{\sim}\THom(\Theta,
  \Phi\circ \Xi)\\
  \pic{AdjunctionIso1}&\mapsto \pic{AdjunctionIso2}
\end{align*}
The inverse isomorphism uses the counit in the obvious way.

In a similar fashion, for $\widehat\Theta\colon Y\to Z$ and $\widehat\Xi\colon
X\to Z$ two 1-morphisms, one obtains an isomorphism
\[\THom(\widehat\Theta\circ\Phi,\widehat\Xi)\xrightarrow{\sim}\THom(\widehat\Theta,
\widehat\Xi\circ\Psi),\]
for which the reader is encouraged to draw the relevant pictures.  It
is worth noting that with respect to the previous isomorphism, $\Psi$
and $\Phi$ have swapped sides in all senses.

Adjunctions are unique up to a canonical isomorphism by the usual
argument.  This means that if $\Psi$ and $\Psi'$ are, say, both left
adjoint to $\Phi$, then there is a canonical isomorphism
$\Psi\stackrel{\sim}\Longrightarrow\Psi'$.  This is pictured below and
it is easy to check that this is an isomorphism.
\[\pic{AdjCanIso}\]

Adjunctions are natural in the sense that they are preserved by
2-func\-tors, so, for instance, given a pair of adjoint kernels in
$\gVar$, the corresponding integral transforms are adjoint functors.

 \subsection{Left and right adjoints of kernels.}
\label{subsec:leftrightadjoints}
\label{subsec:yphimphi}
In an arbitrary 2-category a given 1-morphism might or might not have
a left or a right adjoint, but in the 2-category $\gVar$ every
1-morphism, that is every kernel, has both a left and a right adjoint.
We will see below that for a kernel $\Phi\colon X\to Y$ there are
adjunctions
\[\Phi^\chk\circ\SerK_Y\adjoint\Phi\adjoint \SerK_X\circ
\Phi^\chk,\]
where $\Phi^\chk\colon Y\to X$ means the object $\sHom_{\D(X\times
  Y)}(\Phi,\cO_{X\times Y})$ considered as an object in $\D(Y\times
X)$.  This should be compared with the fact that if $M$ is an
$A$-$B$-bimodule then $M^\chk$ is naturally a $B$-$A$-bimodule.  We
shall see that the two adjunctions above are related in some very
useful ways.

\begin{Proposition}\label{prop:deltalstar}
  If $X$ is a space and $\Delta\colon X\to X\times X$ is the diagonal
  embedding then $\Delta_*\colon \D(X)\to \D(X\times X)$, the
  push-forward on derived categories, is a monoidal functor where
  $\D(X)$ has the usual monoidal tensor product $\otimes$ and
  $\D(X\times X)$ has the composition $\circ$ as the monoidal
  structure.
\end{Proposition}

\begin{Proof}
The proof is just an application of the projection formula.
\end{Proof}

\noindent This has the following immediate consequence.

\begin{Lemma}
  If $\cE$ and $\cF$ are dual as objects in $\D(X)$ then $\Delta_*\cE$
  and $\Delta_*\cF$ are both left and right adjoint to each other as
  1-morphisms in $\gVar$.
\end{Lemma}

In order not to hold-up the flow of the narrative, the proofs of the
remaining results from this section have been relegated to
Appendix~\ref{sec:duality}.

 We begin with some background on the
Serre kernel $\SerK_X$.
Recall from Section~\ref{sec:serre} that the anti-Serre kernel
$\SerK^{-1}_X$ is defined to be $\Delta_*\omega^{-1}_X[-\dim X]$, and
that the Serre kernel $\SerK_X$ is denoted by a dashed-dotted line,
while the anti-Serre kernel $\SerK^{-1}_X$ is denoted by a
dashed-dotted line with a horizontal bar.  The above proposition means
that we actually have maps
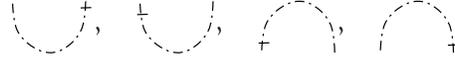
\[\pic{SerK1},\quad
\pic{SerK2},\quad
\pic{SerK3},\quad
\pic{SerK4}
\]
satisfying the following relations
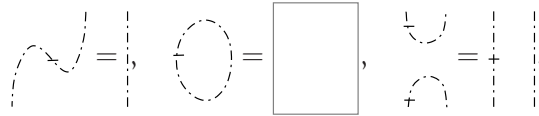
\[\pic{SerK5}=\pic{SerK6},\quad
\pic{SerK7}=\pic{SerK8},\quad
\pic{SerK9}=\pic{SerK10},\quad
\]
and all obvious variations thereof.

In the appendix it is shown that for a kernel $\Phi$ there are natural
morphisms $\rote_\Phi\colon \Phi\circ\SerK_X\circ\Phi^\chk\to \cO_\Delta$
and $\gamma_\Phi\colon \SerK_X^{-1}\to \Phi^\chk\circ \Phi$, denoted in the
following fashion, where the solid, upward oriented lines are labelled
with $\Phi$ and the solid, implicitly downward oriented lines are
labelled by $\Phi^\chk$:
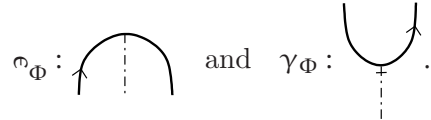
\[\rote_\Phi:\pic{M23}\quad \text{and}\quad \gamma_\Phi:\pic{Y23c}.\]
The main property of $\rote_\Phi$ and $\gamma_\Phi$ is that if we define
$\epsilon_\Phi$, $\overline \epsilon_\Phi$, $\eta_\Phi$ and
$\overline\eta_\Phi$ via
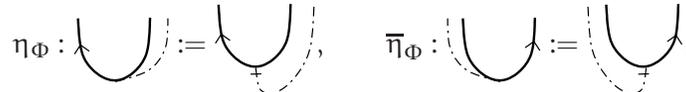
\begin{gather*}
 \eta_\Phi: \pic{eta1}:=\pic{YbentRight},\qquad
\overline\eta_\Phi: \pic{eta2}:=\pic{YbentLeft},\\
\epsilon_\Phi:=\rote_{\Phi^\chk}=\pic{M14},\qquad
\overline\epsilon_\Phi:=\rote_{\Phi}=\pic{M23},
\end{gather*}
then these are the units and counits of adjunctions
\[\Phi^\chk\circ\SerK_Y\adjoint\Phi\adjoint \SerK_X\circ
\Phi^\chk.\]

\subsection{Partial traces.}
We can now define the important notion of partial traces.
\begin{Definition}
   For a kernel $\Phi\colon X\kerto Y$, and
  kernels $\Psi,\Theta\colon Z\kerto X$ define the \emph{left partial trace}
  \[ \THom(\Phi\circ\Theta,\SerK_Y\circ\Phi\circ\Psi) \to
  \THom(\Theta,\SerK_X\circ\Psi) \]
  as
  \[ \pic{UnpartialTraceLeft}\mapsto\pic{PartialTraceLeft}.\]
  Similarly we define a \emph{right partial trace}
  \[ \THom(\Theta'\circ\Phi,\Psi'\circ\Phi\circ\SerK_X) \to
  \THom(\Theta',\Psi'\circ\SerK_Y)\]
  as
  \[ \pic{UnpartialTraceRight}\mapsto\pic{PartialTraceRight}.\]
\end{Definition}
\medskip

\noindent
The following key result, proved in the appendix, says that taking partial
trace does not affect the Serre trace.

\begin{Theorem}
  For a kernel $\Phi\colon X\kerto Y$, a kernel $\Psi\colon Z\kerto X$
  and a kernel morphism
  $\alpha\in\Hom(\Phi\circ\cF,\SerK_Y\circ\Phi\circ\Psi\circ\SerK_Z)$
  then the left partial trace of $\alpha$ has the same Serre trace as
  $\alpha$, i.e., pictorially
  \[\Tr\left(\pic{ApPT1}\right)=\Tr\left(\pic{ApPT2}\right).\]
  The analogous result holds for the right partial trace.
\end{Theorem}

 \subsection{Adjunction as a 2-functor.}
\label{subsec:adj2functors}
As shown in Section~\ref{subsec:leftrightadjoints}, in the
2-cate\-gory $\gVar$ every 1-morphism, that is every kernel,
$\Phi\colon X\to Y$ has a right adjoint $\SerK_X\circ \Phi^\chk$.
This can be extended to a `right adjunction 2-functor' $\tau_R\colon
\gVar^\text{coop}\to \gVar$, where $\gVar^\text{coop}$ means the
contra-opposite 2-category of $\gVar$, which is the 2-category with
the same collections of objects, morphisms and 2-morphisms, but in
which the direction of the morphisms and the 2-morphisms are reversed.

Before defining $\tau_R$, however, it is perhaps useful to think of
the more familiar situation of a one-object 2-category with right
adjoints, i.e., a monoidal category with (right) duals.  So if $\gC$
is a monoidal category in which each object $a$ has a dual $a^\chk$
with evaluation map $\epsilon_a\colon a^\chk\otimes a\to \one$ and
coevaluation map $\eta_a\colon \one \to a\otimes a^\chk$, then for any
morphism $f\colon a\to b$ define $f^\chk\colon b^\chk\to a^\chk$ to be
the composite:
\[b^\chk\xrightarrow{\id\otimes \eta} b^\chk\otimes a\otimes
a^\chk\xrightarrow{\id\otimes f\otimes \id} b^\chk\otimes b\otimes
a^\chk\xrightarrow{\epsilon\otimes \id} a^\chk.\]
This gives rise to a functor $(-)^\chk\colon \gC^\text{op}\to \gC$.

Now return to the case of interest and define $\tau_R\colon
\gVar^\text{coop}\to \gVar$ as follows.  On spaces define
$\tau_R(X):=X$.  On a kernel $\Phi\colon X\to Y$ define
$\tau_R(\Phi):=\SerK_X\circ \Phi^\chk$.  Finally, on morphisms of
kernels define it as illustrated:
\[\tau_R\left(\pic{BeforetauRalpha}\right):=\pic{tauRalpha}.\]
It is a nice exercise for the reader to check that this is a
2-functor.

Clearly a left adjoint 2-functor $\tau_L\colon \gVar^\text{coop}\to
\gVar$ can be similarly created by defining it on a kernel $\Phi\colon
X\to Y$ by $\tau_L(\Phi):=\Phi^\chk\circ\SerK_Y$ and by defining it on
morphisms of kernels by
\[\tau_L\left(\pic{BeforetauRalpha}\right):=\pic{tauLalpha}.\]

\section{Induced maps on homology}
\label{sec:defmaps}
In this section we define $\HH_\blob(X)$, the Hochschild homology of a
space $X$, and show that given a kernel $\Phi\colon X\kerto Y$ we get
pull-back and push-forward maps, $\Phi^*\colon \HH_\blob(Y)\to
\HH_\blob(X)$ and $\Phi_*\colon \HH_\blob(X)\to \HH_\blob(Y)$, such
that if $\Phi$ is right adjoint to $\Psi$ then $\Psi^*=\Phi_*$.

 \subsection{Hochschild cohomology.}
First recall that for a space $X$, one way to define its Hochschild
\emph{cohomology} is as
\[\HH^\blob(X):=\Ext_{X\times X}^\blob(\OO_\Delta,\OO_\Delta).\]
However, the ext-group is just the hom-set $\Hom^\blob_{\D(X\times
  X)}(\OO_\Delta,\OO_\Delta)$ which by the definition of $\gVar$ is
just $\THom^\blob_\gVar(\Id_X, \Id_X)$.  In terms of diagrams, we can
thus denote an element $\phi\in \HH^\blob(X)$ as
\[\pic{phiinHH}.\]
Note that the grading is not indicated in the picture, but this should
not give rise to confusion.

 \subsection{Hochschild homology.}
Now we define $\HH_\blob(X)$ the Hochschild \emph{homology} of a space
$X$ as follows:
\[\HH_\blob(X):=\THom_\gVar^\blob(\SerK_X^{-1},\Id_X) \]
or, in more concrete terms, 
\[ \HH_\blob(X) = \Ext_{X\times X}^{-\blob}(\SerK^{-1}_X,\OO_\Delta).\]
Thus an element $w\in \HH_\blob(X)$ will be denoted
\[\pic{winHH},\]
where again the shifts are understood.

It is worth taking a moment to compare this with other definitions of
Hochschild homology, such as that of Weibel~\cite{Weibel}.  He defines
the Hochschild homology of a space $X$ as
$\H^\blob(X,\Delta^*\cO_\Delta)$, where as usual by $\Delta^*$ we mean
the left-derived functor.  This cohomology group is naturally
identified with the hom-set $\Hom^\blob_{X}(\cO_X,\Delta^*\cO_\Delta)$
which is isomorphic to $\Hom^\blob_{X\times
  X}(\Delta_!\cO_X,\cO_\Delta)$ where $\Delta_!$ is the left-adjoint
of $\Delta^*$.  Direct calculation shows that $\Delta_!\cO_X\cong
\SerK_X^{-1}$ and so our definition is recovered.  Another feasible
definition of Hochschild homology is $\H^\blob(X\times
X,\cO_\Delta\otimes \cO_\Delta)$, and this again is equivalent to our
definition as there is the isomorphism $\SerK_X^{-1}\cong
\cO_\Delta^\chk$.

 \subsection{Push-forward and pull-back.}
\label{def:defmaps}
For spaces $X$ and $Y$ and a kernel $\Phi\colon X\kerto Y$ define the
push-forward on Hochschild homology $\Phi_*\colon \HH_\blob(X)\to \HH_\blob(Y)$ as follows
\[\Phi_*\left(\pic{winHH}\right):=\pic{PhiPushOfw}.\]
For the reader still unhappy with diagrams, for $v\in
\Hom^\blob(\SerK_X^{-1},\Id_X)$, define
$\Phi_*(v)\in \Hom^\blob(\SerK_Y^{-1},\Id_Y)$ as
the following composite, which is read from the above diagram by
reading upwards from the bottom:
\begin{align*}
  \SerK^{-1}_Y
  &\xrightarrow{\gamma}
  \Phi\circ \Phi^\chk
  \xrightarrow{\Id\circ\eta\circ\Id}
  \Phi\circ\SerK_X^{-1}\circ\SerK_X\circ \Phi^\chk
  \xrightarrow{\Id\circ v\circ\Id\circ\Id} \Phi\circ\SerK_X\circ
  \Phi^\chk \xrightarrow{\rote} \Id_Y.
\end{align*}
\medskip

\noindent
Similarly define the pull-back $\Phi^*\colon \HH_\blob(Y)\to
\HH_\blob(X)$ as follows
\[\Phi^*\left(\pic{vinHH}\right):=\pic{PhiPullOfv}.\]
These operations depend only on the isomorphism class of the kernel as
shown by the following.

\begin{Proposition}
\label{prop:indepiso}
  If kernels $\Phi$ and $\hat\Phi$ are isomorphic then they give rise to
  equal push-forwards and equal pull-backs: $\Phi_*=\hat\Phi_*$ and
  $\Phi^*=\hat\Phi^*$.
\end{Proposition}

\begin{Proof}
  This follows immediately from the fact that the 2-morphisms $\gamma_\Phi$,
  $\gamma_{\hat\Phi}$, $\rote_\Phi$ and $\rote_{\hat\Phi}$ of
  Section~\ref{subsec:yphimphi} are natural and thus commute with the
  given kernel isomorphism $\Phi\cong\hat\Phi$.
\end{Proof}
\medskip

\noindent
The push-forward and pull-back operations are functorial in the
following sense.

\begin{Theorem}[Functoriality]
\label{thm:functorial}
  If $\Phi\colon X\to Y$ and $\Psi\colon Y\to Z$ are kernels then the
  push-forwards and pull-backs compose appropriately, namely:
  \[(\Psi\circ\Phi)_*=\Psi_*\circ \Phi_*\colon
  \HH_\blob(X)\to\HH_\blob(Z)\]
  and
  \[(\Psi\circ\Phi)^*=\Phi^*\circ \Psi^*\colon
  \HH_\blob(Z)\to\HH_\blob(X).\]
\end{Theorem}

\begin{Proof}
  This follows from the fact that the right adjunction $\tau_R$ is a
  2-functor.  The adjoint of $\Psi\circ\Phi$ is canonically
  $\tau_R(\Phi)\circ\tau_R(\Psi)$, i.e., is
  $\SerK_X\circ\Phi^\chk\circ\SerK_Y\circ\Psi^\chk$.  This means that
  the unit of the adjunction $\Id\Rightarrow
  \Psi\circ\Phi\circ\SerK_X\circ\Phi^\chk\circ\SerK_Y\circ\Psi^\chk$
  is given by the composition
  $\Id\Rightarrow\Psi\circ\SerK_Y\circ\Psi^\chk\Rightarrow
  \Psi\circ\Phi\circ\SerK_X\circ\Phi^\chk\circ\SerK_Y\circ\Psi^\chk$.
  This gives
  \[(\Psi\circ\Phi)_*(w)=\pic{PsiOfPhiOfw}=\Psi_*\left(\Phi_*(w)\right).
  \qedhere
  \]
\end{Proof}

\begin{Theorem}
\label{thm:adjadj}
  If $\Phi\colon X\to Y$ and $\Psi\colon Y\to X$ are adjoint kernels,
  $\Phi\dashv\Psi$, then we have
  \[\Phi_*=\Psi^*\colon \HH_\blob(X)\to\HH_\blob(Y).\]
\end{Theorem}

\begin{Proof}
  By the uniqueness of adjoints we have a canonical isomorphism,
  $\Psi\cong\tau_R(\Phi)$ and by Proposition~\ref{prop:indepiso}
  $\Psi^*=(\tau_R(\Phi))^*$.  It therefore suffices to show that
  $\Phi_*=(\tau_R(\Phi))^*$.  Observe that
  \[\tau_R(\tau_R(\Phi))=\tau_R(\SerK_X\circ\Phi^\chk)\cong
  \tau_R(\Phi^\chk)\circ
  \tau_R(\SerK_X)\cong\SerK_Y\circ\Phi^{\chk\chk}\circ \SerK_X^{-1},\]
  and similarly
 \[\tau_L(\tau_R(\Phi))\cong
 \Phi\circ\SerK_X^{-1}\circ \SerK_X.\] Of course the latter is
 isomorphic to $\Phi$ but the Serre kernels are left in so to make the
 adjunctions more transparent.  We now get the unit for adjunction
 $\tau_R(\Phi)\dashv\tau_R(\tau_R(\Phi))$ and the counit for the
 adjunction $\tau_L(\tau_R(\Phi))\dashv\tau_R(\Phi)$ as follows:
 \[\pic{tauRunit};\quad\pic{tauRcounit}.\]
 Thus
 \[\tau_R(\Phi)^*\left(w\right)
 =\pic{MessyPhiPushOfw}= \pic{PhiPushOfw}=\Phi_*(w). \qedhere\]
\end{Proof}

\section{The Mukai pairing and adjoint kernels}
\label{sec:adjoint}

In this section we define the Mukai pairing on the Hochschild homology
of a space and show that the push-forwards of adjoint kernels are
themselves adjoint linear maps with respect to this pairing.

First observe from Section~\ref{subsec:adj2functors} that we have two
isomorphisms:
\[\tau_R,\tau_L\colon\HH_\blob(X) = \Hom^{-\blob}(\SerK_X^{-1}, \Id_X) \xrightarrow{\sim}\Hom^\blob(\Id_X, \SerK_X),\] 
given by
\[\tau_R\left(\pic{vinHH}\right):=\pic{tauRv}\qquad\text{and} \quad
\tau_L\left(\pic{vprimeinHH}\right):=\pic{tauLvprime}.\]
Note that this differs slightly from the given definition, but we have used the uniqueness of adjoints.  The above isomorphisms allow the definition of the Mukai pairing as
follows.

\begin{Definition} 
  The \emph{Mukai pairing} on the Hochschild homology of a space $X$ is the
  map
  \[ \MP{{-}}{{-}} \colon \HH_\blob(X) \otimes \HH_\blob(X) \ra \k, \]
  defined by
  \[ \MP{v}{v'}:= \Tr\left(\tau_R(v) \circ \tau_L(v')\right). \]
  Diagrammatically, this is
  \[\MP{\pic{vinHH}}{\pic{vprimeinHH}}:=
  \Tr\left(\pic{tauRv}\pic{tauLvprime}\right).\]
\end{Definition}
\medskip

\noindent
Observe that as $\tau_R$ and $\tau_L$ are both isomorphisms and as the
Serre pairing is non-degenerate, it follows that the Mukai pairing is
non-degenerate.

We can now easily show that adjoint kernels give rise to adjoint maps
between the corresponding Hoch\-schild homology groups.

\begin{Theorem}[Adjointness]
\label{thm:adjointness}
  If $\Phi\colon X\to Y$ and $\Psi\colon Y\to X$ are adjoint kernels,
  $\Psi\dashv\Phi$, then the corresponding push forwards are adjoint
  with respect to the Mukai pairing in the sense that for all $w\in
  \HH_\blob(X)$ and $v\in \HH_\blob(Y)$ we have
  \[\MP{\Psi_*(v)}{w}=\MP{v}{\Phi_*(w)}.\]
\end{Theorem}

\begin{Proof}
  Note first that $\Psi_*=\Phi^*$, by Theorem~\ref{thm:adjadj}.  Thus
  \begin{align*}
    \MP{\Psi_*(v)}{w}&=\MP{\Phi^*(v)}{w}
    =\Tr\left(\pic{tauRPhiv}\pic{tauLw}\right)\\
    &=\Tr\left(\pic{tauRv}\pic{IdPhiarrow}\pic{tauLw}\right)
    =\Tr\left(\pic{tauRv}\pic{tauLPhiw}\right)\\
    &=\MP{v}{\Phi_*(w)}.\qedhere
  \end{align*}
\end{Proof}

\begin{Corollary}
  If the integral kernel $\Phi\colon X \ra Y$ induces an equivalence
  on derived categories, then $\Phi_*\colon\HH_\blob(X) \ra
  \HH_\blob(Y)$ is an isometry.
\end{Corollary}

\begin{Proof}
  If $\Phi$ induces an equivalence, then it has a left adjoint
  $\Psi\colon Y\to X$ which induces the inverse, so
  $\Psi\circ\Phi\cong \Id_X$, and we know that $(\Id_X)_*$
  is the identity map.  Thus
  \[ \MP{\Phi_* v}{\Phi_* w} = \MP{\Psi_* \Phi_* v}{w} =
  \MP{(\Psi\circ \Phi)_* v}{w}= \MP{v}{w}.\qedhere \]
\end{Proof}

\section{The Chern character}
\label{sec:chernchar}

In this section we define the Chern character map $\ch\colon \K_0(X)
\ra \HH_0(X)$.  We discuss the relationship between our construction
and the one of Markarian~\cite[Definition~2]{Mar}.  Then we show that
the Chern character maps the Euler pairing to the Mukai pairing: we
call this the Semi-Hirzebruch-Riemann-Roch Theorem.

\subsection{Definition of the Chern character.}
Suppose $X$ is a space, and $\cE$ is an object in $\D(X)$. Consider
$\cE$ as an object of $\D(\pt\times X)$, i.e., as a kernel $\pt \ra
X$, so there is an induced linear map
\[ \cE_*\colon\HH_\blob(\pt) \ra \HH_\blob(X). \]
Now, because the Serre functor on a point is trivial, $\HH_0(\pt)$ is
canonically identifiable with $\Hom_{\pt\times \pt}(\cO_\pt, \cO_\pt)$
so there is a distinguished class $\bone \in \HH_0(\pt)$ corresponding
to the identity map.  Define the Chern character of $\cE$ as
\[ \ch(\cE) := \cE_*(\bone) \in \HH_0(X). \]
Graphically this has the following description:
\[\ch(\cE):=\pic{chE}.\]
Naturality of push-forward leads to the next theorem.

\begin{Theorem}
\label{thm:commut}
If $X$ and $Y$ are spaces and $\Phi\colon X\to Y$ is a kernel then the
diagram below commutes.
\[
\begin{diagram}[height=2em,width=2em,labelstyle=\scriptstyle]
\D(X) & \rTo^{\Phi\circ {-}} & \D(Y) \\
\dTo_{\ch} & & \dTo_{\ch} \\
\HH_0(X) & \rTo^{\Phi_*} & \HH_0(Y).
\end{diagram}
\]
\end{Theorem}

\begin{Proof}
  Let $\cE$ be an object of $\D(X)$.  We will regard it either as an
  object in $\D(X)$, or as a kernel $\pt \ra X$, and similarly we will
  regard $\Phi\circ\cE$ either as an object in $\D(Y)$ or as a kernel
  $\pt \ra Y$.  By Theorem~\ref{thm:functorial} we have
  \begin{align*}
    \Phi_* \ch(\cE) & = \Phi_* \left(\cE_* (\bone)\right) = (\Phi\circ
    \cE)_* \bone = \ch(\Phi\circ\cE).\qedhere
  \end{align*}
\end{Proof}

 \subsection{The Chern character as a map on K-theory.}
\label{subsec:iotae}
To show that the Chern character descends to a map on K-theory we give
a characterization of the Chern character similar to that of
Markarian~\cite{Mar}.

For any object $\cE\in \D(X)$, which is to be considered an object of
$\D( \pt\times X)$, there are the following two maps:
\[\iota_\cE\colon \HH_\blob(X)\to \Hom_{\D(X\times X)}^\blob(\cE,\SerK_X\circ \cE);
\qquad\pic{vinHH}\,\mapsto \,\pic{iEofvinHH}\]
and
\[\iota^\cE\colon \Hom_{\D(X\times X)}^\blob(\cE,\cE)\to\HH_\blob(X);
\qquad\pic{einHomEE}\mapsto \,\pic{homologyiEofeinHomEE}. \]

Recall that the Mukai pairing is a non-degenerate pairing on $\HH_\blob(X)$ and that the Serre pairing is a perfect pairing between $\Hom_{\D(X\times X)}^\blob(\cE,\SerK_X\circ \cE)$
and $\Hom_{\D(X\times X)}^\blob(\cE,\cE)$.
 With respect to these pairings the two maps $\iota^\cE$ and
$\iota_\cE$ are adjoint in the following sense.

\begin{Proposition}\label{prop:iotasadjoint}
  For $\phi\in \Hom^\blob(\cE,\cE)$ and $v\in \HH_\blob(X)$ the
  following equality holds:
  \[\langle v,\, \iota^\cE\phi\rangle_{\rm M} =\SP{\iota_\cE v}{\phi}\]
\end{Proposition}

\begin{Proof}
  Here in the third equality we use the invariance of the Serre trace
  under the partial trace map.
  \begin{align*}
    \MP{v}{\iota^\cE\phi}&=
    \MPbig{\pic{vinHH}}{\pic{homologyiEofeinHomEE}}
    :=\Tr\left(\pic{iEadjoint3}\right)\\&=
    \Tr\left(\pic{iEadjoint4}\right)
    =\Tr\left(\pic{iEadjoint5}\right)\\
    &=: \SPbig{\pic{iEofvinHH}\quad}{\quad\pic{einHomEE}} =:\SP{\iota_\cE
      v}{\phi}.\qedhere
  \end{align*}
\end{Proof}
\noindent
Note that, using this, the Chern character could have been defined as
\[\ch(\cE):=\iota^\cE(\id_\cE).\]
Then from the above proposition the following is immediate.

\begin{Lemma}
  For $v\in \HH_0(X)$ and $\cE\in \D(X)$ there is the equality
  \[\MP{v}{\ch(\cE)}=\Tr(\iota_\cE(v)),\]
  and this defines $\ch(\cE)$ uniquely.
\end{Lemma}

The fact that the Chern character descends to a function on the
K-group can now be demonstrated.

\begin{Proposition}
  For $\cE\in \D(X)$ the Chern character $\ch(\cE)$ depends only on
  the class of $\cE$ in $\K_0(X)$.  Thus the Chern character can be
  considered as a map
  \[ \ch\colon \K_0(X) \ra \HH_0(X). \]
\end{Proposition}

\begin{Proof}
  It suffices to show that if $\cF\ra \cG\ra \cH\ra \cF[1]$ is an
  exact triangle in $\D(X)$, then
  \[ \ch(\cF) - \ch(\cG) + \ch(\cH) = 0 \]
  in $\HH_0(X)$.

  For $\cG,\cH\in \D(X)$, $\alpha \colon \cG\to \cH$ and $v\in
  \HH_\blob(X)$ the diagram on the left commutes as it expresses the
  equality on the right:
  \[
  \begin{diagram}[height=2em,width=2em]
    \cG & \rTo^{\alpha} & \cH\\
    \dTo^{\iota_\cG( v)} & &
    \dTo_{\iota_\cH( v)}  \\
    \SerK_X \circ\cG & \rTo^{\Id_{\SerK_X}\circ \alpha} & \SerK_X
    \circ\cH.
  \end{diagram}\qquad;\qquad\pic{ChernOnK1}\ =\ \pic{ChernOnK2}\ .
  \]
  In other words, from an element $v\in \HH_\blob(X)$ we get
  $\tau_R(v)\in \Hom(\Id_X,\SerK_X)$, which in turn gives rise to
  a natural transformation between the functors $\Id_{\D(X)}\colon \D(X)\to
  \D(X)$ and $\SerK_X\circ{-}\colon \D(X)\to \D(X)$.  This leads to a
  map of triangles
  \[
  \begin{diagram}[height=2em,width=2em]
    \cF & \rTo & \cG & \rTo & \cH & \rTo & \cF[1] \\
    \dTo_{\iota_\cF( v)} & & \dTo_{\iota_\cG( v)} & &
    \dTo_{\iota_\cH( v)} & & \dTo_{\iota_\cF( v)[1]} \\
    S_X \cF & \rTo & S_X \cG & \rTo & S_X \cH & \rTo & S_X \cF[1].
  \end{diagram}
  \]
  Observe that if we represent the morphism $v$ by an actual map of
  complexes of injectives, and the objects $\cF$, $\cG$ and $\cH$ by
  complexes of locally free sheaves, then the resulting maps in the
  above diagram commute on the nose (no further injective or locally
  free resolutions are needed), so we can apply
  \cite[Theorem~1.9]{May} to get
  \[ \Tr_X(\iota_\cF( v)) - \Tr_X(\iota_\cG( v)) + \Tr_X(\iota_\cH(
  v)) = 0. \] Therefore, by the lemma above, for any $ v\in
  \HH_\blob(X)$,
  \[ \MP{v}{\ch(\cF)-\ch(\cG)+\ch(\cH)}= 0 .\]
  Since the Mukai pairing on $\HH_0(X)$ is
  non-degenerate, we conclude that
  \[ \ch(\cF) - \ch(\cG) + \ch(\cH) = 0. \qedhere\]
\end{Proof}

 \subsection{The Chern character and inner products.}
One reading of the Hirzebruch-Riemann-Roch Theorem is that it says
that the usual Chern character map $\ch\colon \K_0\to \H^\blob(X)$ is
a map of inner product spaces when $\K_0(X)$ is equipped with the
Euler pairing (see below) and $\H^\blob(X)$ is equipped with the
pairing $\langle x_1,x_2\rangle:=(x_1\cup x_2 \cup \Td_X)\cap [X]$.
It is shown in~\cite{CalHH2} that the Hochschild homology Chern
character composed with the Hochschild-Kostant-Rosenberg map
$I_{\text{HKR}}$ gives the usual Chern character:
\[\K_0\xrightarrow{\ch}\HH_0(X)\xrightarrow{I_\text{HKR}}
\bigoplus_p \H^{p,p}(X).\]
Here we show that the Hochschild homology Chern character is an
inner-product map when $\HH_\blob(X)$ is equipped with the Mukai
pairing.

First recall that the Euler pairing on $\K_0(X)$ is defined by
\[\chi(\cE,\cF):=\sum_i (-1)^i \dim \Ext^i_X(\cE, \cF). \]

\begin{Theorem}[Semi-Hirzebruch-Riemann-Roch]
\label{thm:hrr}
The Chern character
\[ \ch\colon \K_0\to \HH_0(X) \]
is a map of inner-product spaces: in other words, for $\cE, \cF\in
\D(X)$ we have
\[ \MP{\ch(\cE)}{\ch(\cF)} = \chi(\cE, \cF). \]
\end{Theorem}

\begin{Proof}
  The first thing to do is to get an interpretation of the Euler
  pairing.  Considering the totality of $\Ext$-groups
  $\Ext^\blob(\cE,\cF)$ as a graded vector space, the Euler
  characteristic is just the graded dimension of
  $\Ext^\blob(\cE,\cF)$, which is to say it is the trace of the
  identity map on $\Ext^\blob(\cE,\cF)$.  Moreover, if ${\mathbf R}\pi_*\colon
  \D(X)\to \D(\pt)$ is the derived functor coming from the map $X\to
  \pt$, then
  $$\Ext^\blob(\cE,\cF)\cong
  \H^\blob(X,\cE^\chk\otimes\cF)\cong {\mathbf R}\pi_*(\cE^\chk\otimes\cF)$$ and
  the latter is just the composition $\cE^\chk\circ \cF$, where
  $\cE^\chk$ and $\cF$ are considered as kernels respectively $X\to
  \pt$ and $\pt\to X$.  Thus, using the invariance of the Serre trace
  under the partial trace map,
  \begin{align*}
    \chi(\cE,\cF)&=\Tr(\Id_{\cE^\chk\circ \cF})
    =\Tr\left(\pic{HRR1}\right)
    =\Tr\left(\pic{HRR2}\right)\\&=\MPbig{\pic{chE}}{\pic{chF}}
    =\MP{\ch(\cE)}{\ch(\cF)}.\qedhere
  \end{align*}
\end{Proof}

 \subsection{Example.}  To have a non-commutative example at
hand, consider the case when $G$ is a finite group acting trivially
on a point.  The orbifold $BG$ is defined to be the global quotient $[\scdot/G]$ and then the category of coherent sheaves on the
orbifold $BG$ is precisely the category of finite
dimensional representations of $G$.  One can naturally identify
$\HH_0(BG)$ with the space of conjugation invariant functions on $G$,
and the Chern character of a representation $\rho$ is precisely the
representation-theoretic character of $\rho$.
See~\cite{WilGr} for details.

\section{Open-closed TQFTs and the Cardy Condition}
\label{sec:cardy}

We conclude with a discussion of open-closed topological field theories in the B-model and we prove that a condition holds for
Hochschild homology which is equivalent to the Cardy Condition in the
Calabi-Yau case.  Appropriate references for open-closed 2d
topological field theories include Moore-Segal~\cite{MooSeg},
Costello~\cite{Cos} and Lauda-Pfeiffer~\cite{LaudaPfeifer}.

 \subsection{Open-closed 2d TQFTs.}
Consider the open and closed 2-cobordism category $\TwoCob$ whose
objects are oriented, compact one-manifolds --- in other words,
disjoint unions of circles and intervals --- and whose morphisms are
(diffeomorphism classes of) cobordisms-with-corners between the source
and target one-manifolds.  A morphism can be drawn as a vertical
cobordism, from the source at the bottom to the target at the top.  As
well as parts of the boundary being at the top and the bottom, there
will be parts of the boundary in between, corresponding to the fact
that this is a cobordism with corners.  An example is shown below.
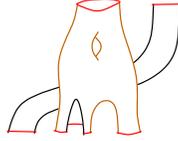
\[\pic{cobordism}\]
Disjoint union makes $\TwoCob$ into a symmetric monoidal category and
an open-closed two-dimensional topological quantum field theory (2d
TQFT) is defined to be a symmetric monoidal functor from $\TwoCob$ to
some appropriate symmetric monoidal target category, which we will
take to be the category of vector spaces or the category of graded
vector spaces.  The category $\TwoCob$ has a simple description in
terms of generators and relations which means that there is a
reasonably straight forward classification of open-closed 2d TQFTs up
to equivalence.  This is what we will now describe.  The following
morphisms generate $\TwoCob$ as a symmetric monoidal category.
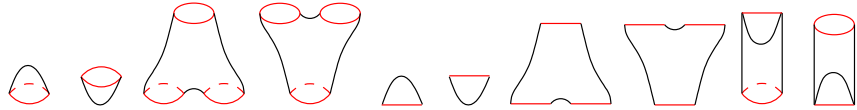
\[\pic{fig1}\]
We will come back to the relations below.

To specify an open-closed 2d TQFT up to equivalence on objects it
suffices to specify the image $\closed$ of the circle and the image
$\open$ of the interval.  The former is called the space of closed
string states and the latter is called the space of open-string
states.  Using the four \emph{planar} generating morphisms pictured
above, together with the relations between them, it transpires that
$\open$, the space of open-string states is precisely a symmetric, but
not-necessarily commutative Frobenius algebra. This means that it is a
unital algebra with a non-degenerate, symmetric, invariant
inner-product.  It is useful to note here that the inner product is
symmetric because the two surfaces pictured below are diffeomorphic,
however these surfaces are \emph{not} ambient isotopic --- so one
cannot be deformed to the other in three-space while the bottom
boundary is fixed.
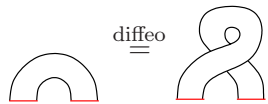
\[\pic{syminnerprod1}\stackrel{\text{diffeo}}{=}\pic{syminnerprod2}\]
On the other hand, the first four generating morphisms, along with
their relations, mean that $\closed$, the space of closed string
states, is a \emph{commutative} Frobenius algebra. The last two
morphisms mean that there are maps $i_*\colon \closed\to\open$ and
$i^*\colon \open \to \closed$, and by the relations these are adjoint
with respect to the pairings on these spaces.  Moreover, $i_*$ is an
algebra map, such that its image lies in the centre of $\open$. The
final relation that these must satisfy is the \emph{Cardy Condition}.
In terms of the generators pictured above this is the following
relation:
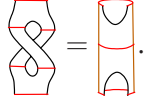
\[\pic{cardy1}=\pic{cardy2}.\]
Note again that these surfaces are diffeomorphic but not isotopic in
three-space.  In terms of maps, writing $\mu\colon \open\otimes
\open\to \open$ and $\delta\colon\open\to\open\otimes\open$ for the
product and coproduct of the open string state space and writing
$\tau\colon\open\otimes\open\to\open\otimes\open$ for the symmetry in
the target category, the Cardy Condition is the equality of maps from
$\open$ to $\open$:
\[\mu\circ \tau\circ \delta=i_*\circ i^*.\]
We will have reason to use an equivalent condition below.

To summarize, having an open-closed 2d TQFT is equivalent to having
the data of a commutative Frobenius algebra $\closed$, a symmetric
Frobenius algebra $\open$, and an algebra map $i_*\colon \closed \to
\open$ with central image, such that the Cardy Condition is satisfied.

 \subsection{Open-closed 2d TQFTs with D-branes.}
A more interesting model of string theory is obtained when we specify
a set of `boundary conditions' or `D-branes' for the open strings. For
a mathematician this just means a set of labels for the boundary
points of objects.  So fix a set $\Lambda$ of labels, and consider the
category $\TwoCobL$ of open-closed cobordisms such that the objects
are compact, oriented one-manifolds with the boundary points labelled
with elements of the set $\Lambda$, and morphisms having their
internal boundaries labelled compatibly with their boundaries. Here is
an example of a morphism from the union of the circle and the interval
labelled $(B,A)$, to the interval labelled $(B,A)$.
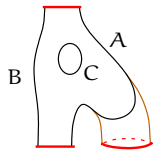
\[\pic{labelledcobordism}\]

Now a $\Lambda$-labelled open-closed TQFT is a symmetric monoidal
functor to some appropriate target category which we will again take
to be the category of vector spaces or the category of graded vector
spaces. Moreover, the category $\TwoCobL$ is similarly generated by
morphisms as listed above, but now they must all be labelled, and the
relations are just labelled versions of the previous relations.  Thus
we can similarly classify $\Lambda$-labelled open-closed TQFTs.  Once
again the image of the circle is a commutative Frobenius algebra,
$\closed$.  However, rather than getting a single vector space $\open$
associated to an interval, we get a vector space $\open_{BA}$
associated to each ordered pair $(B,A)$ of elements of $\Lambda$; so
we do not get a single Frobenius algebra, but rather something which
could be called a `Frobenius algebra with many objects' or a
`Frobenius algebroid', but, for the reason explained below, such a
thing is commonly known as a Calabi-Yau category.  It is a category in
the following sense.  We take the category whose objects are
parametrized by $\Lambda$ and, for $A,B\in \Lambda$, the morphism set
$\Hom(A,B)$ is taken to be $\open_{BA}$ (this is consistent with us
reading diagrams from right to left).  The composition
$\mu_{CBA}\colon\open_{CB}\otimes\open_{BA}\to\open_{CA}$ is given by
the image of the appropriately labelled version of the morphism
pictured.
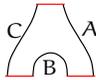
\[\pic{labelledproduct}\]
The Frobenius or Calabi-Yau part of the structure is a --- possibly
graded --- perfect pairing $\open_{AB}\otimes\open_{BA}\to \k$: the
grading degree of this map is called the dimension of the Calabi-Yau
category.

So to specify a labelled open-closed 2d TQFT it suffices to specify a
commutative Frobenius algebra $\closed$, a Calabi-Yau category $\open$
and an algebra map $i_A\colon \closed\to \open_{AA}$ with central
image, for each object $A$, such that the labelled version of the
Cardy Condition holds.

 \subsection{The open-closed 2d TQFT from a Calabi-Yau
  manifold.}
Associated to a Calabi-Yau manifold $X$ there are two standard 2d TQFTs
coming from string theory, imaginatively named the 
A-model and the B-model:  it is the B-model we will be
interested in here.  In the B-model the boundary conditions are
supposed to be ``generated'' by complex submanifolds of $X$ so the boundary conditions are
taken to be complexes of coherent sheaves on $X$;  the open string
category is then supposed to be the derived category of coherent sheaves on
$X$.  This is indeed a Calabi-Yau category, which is why such
categories are so named:  for each $\cE$ and $\cF$, the requisite pairing
$\Hom^\blob_{\D(X)}(\cE,\cF)\otimes \Hom^\blob_{\D(X)}(\cF,\cE)\to
\k[-\dim X]$ comes from the Serre pairing as a Calabi-Yau
manifold is precisely a manifold with a trivial canonical bundle.

According to the physics, the closed string state space $\closed$ should be $\Hom^\blob_{\D(X\times X)}(\cO_\Delta,\cO_\Delta)$, in other words, the
Hochschild cohomology algebra
 $\HH^\blob(X)$.
As $X$ is Calabi-Yau, a trivialization of the canonical bundle induces
an isomorphism between Hochschild cohomology and Hochschild homology,
up to a shift. This means that the closed string space $\closed$
has both the cohomological product and the Mukai pairing, and these make  $\closed$ into a
Frobenius algebra.

We need to specify the algebra maps $i_\cE\colon \closed \to
\open_{\cE\cE}$.  These are maps
\[ i_\cE\colon \Hom^\blob(\cO_\Delta,\cO_\Delta)\to
\Hom^\blob(\cE,\cE) \]
which can be given by interpreting $\cE$  as a kernel $\pt\to X$ and taking $i_\cE$ to be convolution with the identity on $\cE$. This is given
diagrammatically on an element
$\phi\in\Hom^\blob(\cO_\Delta,\cO_\Delta)$ as follows.
\[\pic{phiinHH}\quad\mapsto\quad \pic{iEofphiinHH}.\]

At this point it should be noted that $\closed$ is to be thought of as the centre of the category $\open$.  The notion of centre is generalized from algebras to categories by taking the centre of a category to be the natural transformations of the identity functor; however, in a 2-category an appropriate notion of the centre of an object is the set of 2-endomorphisms of the identity morphism on that object.  This means that $\closed$ is the centre, in this sense, of the category $\open$ in the 2-category $\gVar$.

The map going the other way, $i^\cE\colon \Hom^\blob(\cE,\cE)\to
\Hom^\blob(\cO_\Delta,\cO_\Delta)$ is given by taking the trace, namely for
$e\in \Hom_{\D(X)}(\cE,\cE)$ the map is given by
\[\pic{ereallyinHomEE}\quad\mapsto \quad\pic{iEofeinHomEE}.\]
This definition relies on the fact that $X$ is Calabi-Yau, so that the
Serre kernel is, up to a shift, just the identity 1-morphism $\Id_X$.

%
%

An argument similar to Proposition~\ref{prop:iotasadjoint} shows that
$i_\cE$ and $i^\cE$ are adjoint.  In order to argue that we indeed
have an open-closed TQFT it remains to show that the Cardy Condition
holds.  In fact, we will prove a more general statement, the Baggy
Cardy Condition.

 \subsection{The Baggy Cardy Condition.}
In the case of a manifold $X$ that is not necessarily Calabi-Yau we
don't have the same coincidence of structure as above: we no longer
have a Frobenius algebra $\HH^\blob(X)$; rather we have an algebra
$\HH^\blob(X)$ and an inner product space $\HH_\blob(X)$.  This means
that we can not formulate the Cardy Condition as it stands.  We now state a
condition which makes sense for an arbitrary, non-Calabi-Yau manifold and
which is equivalent to the Cardy Condition in the  Calabi-Yau case.

\begin{Theorem}
  Suppose that $\open$ is a Calabi-Yau category and $\closed$ is an
  inner product space, such that for each $A\in \open$ there are
  adjoint maps $i^A\colon \open_{AA}\to \closed$ and $i_A\colon
  \closed \to \open_{AA}$.  Then the Cardy Condition
  \[\mu_{BAB}\circ\tau\circ\delta_{ABA}=i_B\circ i^A\]
  is equivalent to the following equality holding for all $a\in
  \open_{AA}$ and $b\in \open_{BB}$, where the map ${}_{a}m_{b}\colon
  \open_{AB}\to \open_{AB}$ is the map obtained by pre-composing with
  $a$ and post-composing with $b$:
  \[\left\langle i^B{-},i^A{-}\right\rangle_\closed=
  \Tr{}_{-}m_{-}.\]
\end{Theorem}

\begin{Proof}
  The first thing to do is examine the left-hand side of the Cardy
  Condition.  As $\open$ is a Calabi-Yau category there is the
  following equality of morphisms $\open_{AA}\to \open_{BB}$.
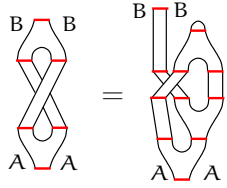
  \[\pic{cardyA}=\pic{cardyB}\]
  Note that this does not require any reference to $\closed$, but it
  does fundamentally require the symmetry of the inner product.  This
  is reflected in the fact that the surfaces underlying the above
  pictures are diffeomorphic but not ambient isotopic.

  This means that the Cardy Condition is equivalent to the following
  equality.
  \[\pic{cardyC}=\pic{cardyD}\]
  By the non-degeneracy of the inner product on $\open_{BB}$ this is
  equivalent to the equality of two maps $\open_{BB}\otimes
  \open_{AA}\to\k$ which are drawn as follows.
  \[\pic{cardyE}=\pic{cardyF}\]
  The right-hand side is instantly identifiable as $\left\langle
    i^B{-},i^A{-}\right\rangle_\closed$.  The left-hand side is
  identifiable as the trace of the triple composition map
  $\open_{BB}\otimes\open_{BA}\otimes\open_{AA}\to \open_{BA}$ which
  gives the required result.
\end{Proof}

We can now show that the alternative condition given in the above
theorem holds for the derived category and Hochschild homology of
\emph{any} space: in particular, the Cardy Condition holds for
Calabi-Yau spaces.

\begin{Theorem}[The Baggy Cardy Condition]
\label{thm:cardy}
Let $X$ be a space, let $\cE$ and $\cF$ be objects in $\D(X)$ and
consider morphisms
\[ e\in \Hom_{\D(X)}(\cE, \cE) \quad\text{and}\quad f\in
\Hom_{\D(X)}(\cF, \cF). \]
Define the operator
\[ {}_f m_e \colon\Hom^\blob_{\D(X)}(\cE, \cF) \ra
\Hom^\blob_{\D(X)}(\cE, \cF) \]
to be post-composition by $f$ and pre-composition by $e$.  Then we
have
\[ \Tr {}_f m_e=\MP{\iota^\cE(e)}{\iota^\cF(f)}, \]
where $\iota^\cE$, $\iota^\cF$ are the maps defined in
Section~\ref{subsec:iotae}, and $\Tr$ denotes the (super) trace.
\end{Theorem}

\begin{Proof}
  The proof is very similar to the proof of the
  Semi-Hirzebruch-Rie\-mann-Roch Theorem (Theorem~\ref{thm:hrr}).  The first thing
  to observe is that $\Hom(\cE,\cF)\cong \cE^\chk\circ \cF$ and that
  ${}_f m_e$ is just $e^\chk\circ f$. However, we have
  $\tau_R(e)=\SerK_\pt\circ e^\chk$ and $\SerK_\pt$ is trivial so
  $e^\chk=\tau_R(e)$.  Putting this together with the invariance of
  the Serre trace under the partial trace we get the following
  sequence, and hence the required result.
  \begin{align*}
    \Tr {}_f m_e&=\Tr\left(\pic{baggy1}\right)
    =\Tr\left( \pic{baggy2}\right)\\
    &=\Tr\left( \pic{baggy3}\right)
    =\Tr\left( \pic{baggy4}\right)\\
    &=\MPbig{\pic{baggy5a}}{\pic{baggy5b}}
    =\MP{\iota^\cE\!(e)}{\iota^\cF\!(f)}.\qedhere
  \end{align*}
\end{Proof}
\medskip

\noindent
Observe that the Semi-Hirzebruch-Riemann-Roch Theorem is a direct
consequence of the Baggy Cardy Condition, with $e=\id_\cE$,
$f=\id_\cF$.

\appendix

\section{Duality and partial trace}
\label{sec:duality}

In this appendix we show that given a kernel $\Phi\colon X\ra Y$ and
its dual kernel \(\Phi^\chk\colon Y \ra X\) there are canonical
2-morphisms
\[ \SerK_X^{-1} \ra \Phi^\chk \circ \Phi\quad\text{and}\quad \Phi
\circ \SerK_X \circ \Phi^\chk \ra \Id_Y \]
giving rise to a variety of natural adjunctions satisfying
a number of compatibility relations.

The notion of duality in $\gVar$ is seen to be a middle-ground between
the operations $\tau_L$ and $\tau_R$: it is an involution, unlike
$\tau_L$ and $\tau_R$, but it does not respect composition, which
$\tau_L$ and $\tau_R$ do.

 \subsection{Polite duality.}
Recall from Section~\ref{subsec:serXY} that for every space $X$ there
is the Serre kernel $\SerK_X\colon X \ra X$ such that for spaces $X$
and $Y$ the functor
\[ \SerK_Y \circ {-} \circ \SerK_X\colon \bHom_\gVar(X, Y) \ra
\bHom_\gVar(X, Y) \]
is a Serre functor for the category $\bHom_\gVar(X,Y)$.

\begin{Definition}
\label{def:dualpair}
If $\Phi\colon X \ra Y$ and $ \pd{\Phi}\colon Y \ra X$ are kernels then
a \emph{polite duality} between them, denoted \( \Phi \dual \pd{\Phi}
\), consists of adjunctions as follows (numbered as shown)
\begin{gather*}
  \pd{\Phi} \circ \SerK_Y \adjoint_1 \Phi \adjoint_2 \SerK_X \circ
  \pd{\Phi},
  \\
  \Phi \circ \SerK_X \adjoint_3 \pd{\Phi} \adjoint_4 \SerK_Y \circ \Phi,
\end{gather*}
such that the following compatibility relations hold.
\begin{quote}
\begin{itemize}
\item[(1+2).] For kernels $\Theta\colon Z\ra X$ and $\Psi\colon Z\ra
  Y$, the diagram of isomorphisms below commutes:
  \[
  \begin{diagram}[height=2em,width=2em,labelstyle=\scriptstyle]
    \Hom(\Psi, \Phi \circ\Theta) & \rTeXto^1 &
    \Hom(\pd{\Phi} \circ\SerK_Y\circ \Psi, \Theta) \\
    \dTeXto^{\text{Serre}} & &
    \dTeXto^{\text{Serre}} \\
    \Hom(\Phi\circ\Theta, \SerK_Y\circ \Psi \circ \SerK_Z)^\chk &
    \rTeXto^2 & \Hom(\Theta, \SerK_X\circ \pd{\Phi}\circ \SerK_Y\circ
    \Psi\circ \SerK_Z)^\chk.
  \end{diagram}
  \]
\item[(2+3).] The composite map
  \[ \Hom(\Theta, \pd{\Phi}\circ\Psi) \xrightarrow[3]{\sim}
  \Hom(\Phi\circ\SerK_X\circ\Theta, \Psi) \xrightarrow[2]{\sim}
  \Hom(\SerK_X\circ\Theta, \SerK_X\circ \pd{\Phi}\circ \Psi) \]
  is the one induced by composition with $\SerK_X$.
\item[(3+4).] Same as (1+2), but for adjunctions 3 and 4.
\item[(1+4).] Same as (2+3), but for adjunctions 1 and 4.
\end{itemize}
\end{quote}
\end{Definition}

It is useful to think of this definition in terms of string diagrams.
For each of the four adjunctions we get a unit and a counit.  Denoting
$\Phi$ by an upward oriented line and $\pd{\Phi}$ by a downward oriented
line, we can draw the units and counits as follows.
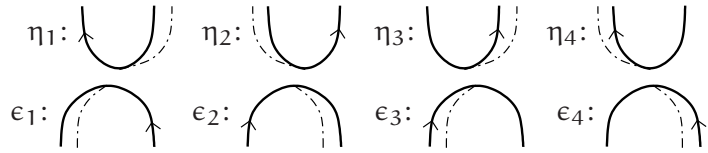
\begin{gather*}
\eta_1\colon\pic{eta1}\quad
\eta_2\colon\pic{eta2}\quad
\eta_3\colon\pic{eta3}\quad
\eta_4\colon\pic{eta4}\\
\epsilon_1\colon\pic{epsilon1}\quad
\epsilon_2\colon\pic{epsilon2}\quad
\epsilon_3\colon\pic{epsilon3}\quad
\epsilon_4\colon\pic{epsilon4}
\end{gather*}

Relation (2+3) can be represented graphically in the following way:
for any $\alpha\in\Hom(\Theta,\pd{\Phi}\circ\Psi)$ we have
\[\pic{R2and3a}=\pic{R2and3b}.\]
From this we can deduce the equality
\[\pic{Y23a}=\pic{Y23b},\]
so we define the following diagram to be this common morphism:
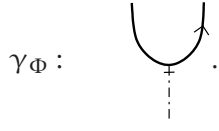
\[\gamma_\Phi:\qquad\pic{Y23c}.\]
Similarly we can deduce that the counits $\epsilon_2$ and $\epsilon_3$
are equal, and we define the following diagram to be the common
morphism:
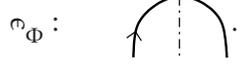
\[\rote_\Phi:\qquad\pic{M23}.\]
Likewise, using the relation (1+4), from the units $\eta_1$ and
$\eta_4$, and the counits $\epsilon_1$ and $\epsilon_4$ we obtain the
common morphisms denoted as follows.
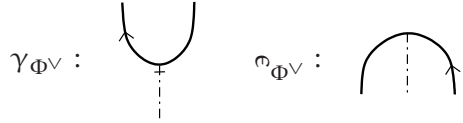
\[\gamma_{\Phi^\chk}:\quad\pic{Y14}\qquad \rote_{\Phi^\chk}:\quad \pic{M14}\]
So the four units and four counits are obtained from these two
M-shaped and two Y-shaped morphisms.

Relations (1+2) and (3+4) in the definition of polite duality are
essentially equivalent to the invariance of the Serre trace under a
partial trace in the following sense.  Given a polite duality
$\Phi\dual\pd{\Phi}$ and a morphism $\alpha\in
\Hom(\Phi\circ\Theta,\SerK\circ \Phi\circ \Psi)$ we can define the
left partial trace in $\Hom(\Theta,\SerK\circ \Psi)$ with respect to
$\pd{\Phi}$ as drawn below.
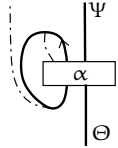
\[\pic{PartialTraceLeft}\]
Similarly we can define a right partial trace when the $\Phi$ is on
the right rather than on the left.  So the following diagram is the
right partial trace with respect to $\Phi$ and $\pd{\Phi}$ of a morphism
$\alpha'\in \Hom(\Theta'\circ\Phi,\Psi'\circ \Phi\circ\SerK)$.
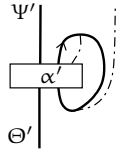
\[\pic{PartialTraceRight}\]

\begin{Theorem}
  If $\Phi$ and $\pd{\Phi}$ form a politely dual pair of kernels, then
  the Serre trace is invariant under partial trace with respect to
  $\Phi$:
  \[\Tr\left(\pic{ApPT1}\right)=\Tr\left(\pic{ApPT2}\right);
  \quad \Tr\left(\pic{ApPT1b}\right)=\Tr\left(\pic{ApPT2b}\right).\]
\end{Theorem}

\begin{Proof}
  We will just prove the left partial trace case.  Relation (1+2) says
  that for $\beta\in \Hom (\Psi,\Phi\circ\Theta)$ and $\gamma\in
  (\Theta,\SerK\circ\pd{\Phi}\circ\SerK\circ\Psi\circ\SerK)$ we have the
  equality
  \[\Tr\left(\pic{ApPT3}\right)=\Tr\left(\pic{ApPT4}\right).\]
  Applying this with
  \[\beta:=\pic{ApPT6}\quad\text{and}\quad
  \gamma:=\pic{ApPT5}\] we obtain the equality
  \[\Tr\left(\pic{ApPT7}\right)=\Tr\left(\pic{ApPT8}\right).\]
  The result follows from the commutativity property of the Serre
  trace, namely that $\Tr(\beta\circv \alpha)=\Tr(\SerF(\alpha)\circv
  \beta)$.
\end{Proof}

\begin{Proposition}
\label{prop:compduals}
  Polite dualities $\Phi \dual \pd{\Phi}$ and $\Psi \dual \pd{\Psi}$ for
  two kernels $\Phi\colon X\ra Y$ and $\Psi\colon Y\ra Z$ canonically
  induce a polite duality
  \[ \Psi \circ \Phi \dual \pd{\Phi} \circ \SerK_Y \circ \pd{\Psi}. \]
\end{Proposition}

\begin{Proof}
  The required adjunctions are constructed in the obvious fashion,
  using the fact that the composition, in reverse order, of adjoints
  (left or right) of composable functors is naturally an adjoint (in
  the same direction) of the composition of the functors.  The
  compatibilities required by the polite duality follow from straight
  forward checks.
\end{Proof}

Adjunctions are intrinsic parts of a polite duality.  Conversely, the
following proposition shows that any adjunction induces a polite
duality in a natural way.

\begin{Proposition}
\label{prop:adjgivesdual}
Given an adjunction $\Psi \adjoint \Phi$, define $\pd{\Phi}:= \Psi \circ
\SerK_Y^{-1}$.  Then there exists a polite duality $\Phi \dual
\pd{\Phi}$ where the adjunction $\adjoint_1$ is the given one
$\Psi\adjoint \Phi$.
\end{Proposition}

\begin{Proof}
  As $\SerK_X\circ \pd{\Phi}= \SerK_X\circ\Psi\circ\SerK^{-1}_Y$, compatibility relation (1+2) yields an adjunction
  \[ \Phi \adjoint_2 \SerK_X
  \circ \pd{\Phi}. \] Similarly, relations (1+4) and (2+3) force
  adjunctions
  \[ \Phi \circ \SerK_X \adjoint_3 \pd{\Phi} \adjoint_4 \SerK_Y \circ
  \Phi, \]
  and it is easy to see that relation (3+4) is then automatically
  satisfied.
\end{Proof}

\begin{Proposition}
\label{prop:existdual}
Let $\Phi\in\D(X\times Y)$ be a kernel, and let
$\Phi^\chk\in\D(Y\times X)$ be the dual $\mathbf{R}\sHom_{X\times
  Y}(\Phi, \cO_{X\times Y})$, regarded as a kernel from $Y$ to $X$.
Then there exists a canonical polite duality \(\Phi \dual \Phi^\chk\).
\end{Proposition}

\begin{Proof}
For kernels $\Theta \in \D(Z\times X)$ and $\Psi\in \D(Z\times Y)$
consider the sequence of isomorphisms
\begin{align*}
\Hom_{\D(Z\times Y)} & (\Psi, \Phi \circ \Theta) \\*
  & \iso \Hom_{\D(Z\times Y)}(\Psi, \pi_{ZY, *}
      (\pi_{ZX}^* \Theta \otimes \pi_{XY}^* \Phi)) \\
  & \iso \Hom_{\D(Z\times X)}(\pi_{ZX,!}
    (\pi_{ZY}^* \Psi \otimes \pi_{YX}^* \Phi^\chk), \Theta) \\
  & = \Hom_{\D(Z\times X)}(\pi_{ZX, *}
      (\pi_{ZY}^*\Psi \otimes \pi_{YX}^*
        (\Phi^\chk \otimes \omega_Y[\dim Y])), \Theta) \\
  & = \Hom_{\D(Z\times Y)}((\Phi^\chk\circ \SerK_Y)\circ\Psi, \Theta).
\end{align*}
Taking $Z = Y$, $\Theta = \Phi^\chk \circ \SerK_Y$ and $\Psi=\Id_Y$
yields a morphism {\em of kernels}
\[ \Id_Y \ra \Phi \circ \Phi^\chk\circ \SerK_Y; \]
in a similar fashion we obtain a morphism
\[ \Phi^\chk \circ \SerK_Y \circ \Phi \ra \Id_X. \]
These two morphisms satisfy the identities needed to make $\Phi^\chk
\circ \SerK_Y$ the left adjoint of $\Phi$.
Proposition~\ref{prop:adjgivesdual} gives the result.
\end{Proof}

 \subsection{Reflexively polite kernels.}
We still need to address one more compatibility between the dualities
constructed above.  Given a kernel $\Phi$, the previous proposition
yields a polite duality
\[ \Phi \dual \Phi^\chk. \]
Given any polite duality $\Phi\dual \pd{\Phi}$, we get, symmetrically, a polite
duality $\pd{\Phi}
\dual \Phi$ by switching the adjunctions $1$ and  $3$ and the adjunctions
$2$ and $4$.  Thus there is a natural polite duality
\[ \Phi^\chk \dual \Phi. \]
On the other hand, applying Proposition~\ref{prop:existdual} to the
kernel $\Phi^\chk$ we get a polite duality
\( \Phi^\chk \dual \Phi^{\chk\chk}\)
and then, using the canonical identification $\Phi^{\chk\chk} \iso
\Phi$, we get {\em another} polite duality
\[ \Phi^\chk \dual \Phi.\]
The fundamental question is whether these two polite dualities are the same.

\begin{Definition}
  Let $\Phi\colon X\ra Y$ be a kernel.  We shall say that $\Phi$ is
  \emph{reflexively polite} if the two dualities above are equal.
\end{Definition}

\noindent Immediately we get the following result.

\begin{Proposition}
  If a kernel $\Phi$ is reflexively polite then so is $\Phi^\chk$.
\end{Proposition}

\begin{Proposition}
\label{prop:politecompduals}
If $\Phi\colon X \ra Y$ and $\Psi\colon Y \ra Z$ are reflexively
polite kernels then so is $\Psi \circ \Phi$.
\end{Proposition}

\begin{Proof}
  This follows at once from the fact that the adjunctions defined in
  Proposition~\ref{prop:compduals} are canonical.
\end{Proof}

\begin{Proposition}
\label{prop:dualbox}
Suppose that $\Phi\colon X\ra Y$ and $\Psi\colon X'\ra Y'$ are
reflexively polite kernels then so is $\Phi\boxtimes \Psi\colon X
\times X' \ra Y\times Y'$, the kernel defined by
\[ \Phi\boxtimes \Psi := \pi_{XY}^* \Phi \otimes \pi_{X'Y'}^*\Psi. \]
\end{Proposition}

\begin{Proof}
  This is obvious once one realises that all operations decompose with
  respect to the box product operation, and the canonical bundle of a
  product of spaces is the box product of the canonical bundles of the
  factors.
\end{Proof}

\begin{Proposition}
\label{prop:dualker}
Let $\Delta\colon X \ra X\times X$ denote the diagonal map, let $\cE$
be an object in $\D(X)$, then the kernel $\Delta_*\cE\colon X \ra X$
is reflexively polite.
\end{Proposition}

\begin{Proof}
  The dual $(\Delta_*\cE)^\chk$ of $\Delta_*\cE$ is naturally
  identified with 
  \[ \Delta_* (\cE^\chk \otimes \omega_X^{-1}[-\dim X]). \] The result
  now follows from Proposition~\ref{prop:deltalstar} and the fact that
  the functor given by tensoring with an object of the derived
  category is left and right adjoint to the one obtained by tensoring
  with the dual object.
\end{Proof}

\begin{Proposition}
\label{prop:dualcox}
Let $X$ be a space, then the kernel $\cO_X$ considered as a kernel
$X\ra \pt$ is reflexively polite, and so also is $\cO_X^\chk\colon
\pt\to X$.
\end{Proposition}

\begin{Proof}
  First, observe that checking reflexive duality amounts to checking
  the equality of two quadruples of adjunctions.  In each quadruple,
  any one of the adjunctions determines the remaining three, and thus
  in order to check the equality it suffices to check that one of the
  compatibilities (1+2), (1+3), {etc}.\ holds with one adjunction
  chosen from one duality, and the other one from the second duality.

  Furthermore, adjunctions are completely determined by their
  respective units and counits, and these are determined by considered
  functors {\em between derived categories} of the form
  \[ \Phi \circ {-} \colon \D(Z\times X) =:\bHom_\gVar(Z, X) \ra
  \bHom_\gVar(Z, \pt) := \D(Z), \]
  and analogues for $\Phi^\chk \circ {-}$, for various choices of the
  space $Z$.  Thus, if we argue that the desired equality holds for
  the adjunctions between these induced functors (for arbitrary choice
  of $Z$), we will have argued that $\Phi$ is reflexively polite.

  For ease of notation, write $\Phi$ for the kernel $\cO_X\colon X\to \pt$.  For a given space $Z$, let $\pi_Z\colon X\times Z \ra Z$ denote the
  projection.  The functors $\Phi \circ {-}$ and $\Phi^\chk \circ
  {-}$ are then naturally identified with $\pi_{Z,*}$ and $\pi_Z^*$,
  respectively.  Adjunctions $1$ and $2$ from the polite duality
  $\Phi\dual \Phi^\chk$ correspond to the classical adjunctions
  \begin{align*}
    \pi_Z^* & \adjoint \pi_{Z,*} \adjoint \pi_Z^!, \intertext{while
      the same ones from the polite duality $\Phi^\chk\dual
      \Phi^{\chk\chk}$ correspond to} \pi_{Z,!} & \adjoint \pi_Z^*
    \adjoint \pi_{Z,*}.
  \end{align*}
  Indeed, the standard way (see~\cite[Theorem 4.6]{CalSkim}) to define
  the adjunctions $\pi_{Z,*} \adjoint \pi_Z^!$ and $\pi_{Z,!}\adjoint
  \pi_Z^*$ is to require them to satisfy the analogue of condition
  (1+2) from the definition of a polite duality.

  The condition of polite duality can now can be stated as the
  statement that the composite isomorphism
  \begin{align*}
    \Hom_{X\times Z}(\Theta, & \pi_Z^* \Psi) \iso \Hom_Z(\pi_{Z,!}
    \Theta, \Psi) =
    \Hom_Z\left(\pi_{Z,*}\SerF(\Theta), \Psi\right) \\
    & \iso \Hom_{X\times Z}\left(\SerF(\Theta), \pi_Z^! \Psi\right)
    = \Hom_{X\times Z}\left(\SerF(\Theta), \SerF (\pi_Z^*
      \Psi)\right)
  \end{align*}
  is the one induced by the functor $\SerF({-}) = {-} \otimes
  \pi_X^* \omega_X[\dim X]$, where $\pi_X$ denotes the projection from
  $X\times Z$ to $X$.  This fact corresponds to the fact that the
  diagram marked (!) below commutes {\small
    \[
    \begin{diagram}[height=2em,width=2em,labelstyle=\scriptstyle]
      \Hom(\Theta, \pi^* \Psi) & \rLine^\sim & \Hom(\SerF ^{-1} \pi_*
      \SerF \Theta, \Psi) &
      \rLine^\sim & \Hom(\pi_* \SerF \Theta, \SerF \Psi) & \rLine^\sim & \Hom(\SerF \Theta, \SerF \pi_* \Psi) \\
      \dLine^\sim & \mbox{Serre} & \dLine^\sim & \mbox{(!)} &
      \dLine^\sim &
      \mbox{Serre} & \dLine^\sim \\
      \Hom(\pi^*\Psi, \SerF \Theta)^\chk & \rLine^\sim & \Hom(\Psi,
      \pi_*\SerF \Theta)^\chk & \rEq & \Hom(\Psi, \pi_*\SerF
      \Theta)^\chk & \rLine^\sim & \Hom(\pi^*\Psi, \Theta)^\chk.
    \end{diagram}
    \]}

\end{Proof}

\begin{Proposition}
\label{prop:pulldualpush}
Let $\pi_X\colon X \times Y \ra X$ be the projection, and (abusing
notation) denote by $\pi_{X,*}\colon X\times Y \ra X$ and
$\pi_X^*\colon X\ra X\times Y$ the kernels represented in $\D(X\times
X \times Y)$ by the structure sheaf of the graph of $\pi_X$.  Then
$\pi_{X,*}$ and $\pi_X^*$ are reflexively polite.
\end{Proposition}

\begin{Proof}
  Both kernels are of the form $\cO_{\Delta_X} \boxtimes \cO_Y$, and
  the result then follows from
  Propositions~\ref{prop:dualcox},~\ref{prop:dualker},
  and~\ref{prop:dualbox}.
\end{Proof}

\begin{Theorem}
  Every kernel is reflexively polite.
\end{Theorem}

\begin{Proof}
  With notation as in Proposition~\ref{prop:pulldualpush}, any kernel
  $\Phi\colon X\to Y$ decomposes as
\[ \Phi = \pi_{Y,*} \circ ({-} \otimes \Phi) \circ \pi_X^*. \]
Since the claim of the theorem holds for each individual kernel in the
decomposition (Propositions~\ref{prop:pulldualpush}
and~\ref{prop:dualker}), the result follows from
Proposition~\ref{prop:compduals}.
\end{Proof}

\end{document}